\documentclass[a4paper]{amsart} 
\usepackage{amsmath}
\usepackage{amssymb}
\usepackage{verbatim}
\usepackage{amsthm}
\usepackage{mathrsfs}
\usepackage{graphicx}
\usepackage{grffile}
\usepackage{mathtools}
\usepackage[colorlinks,pdfdisplaydoctitle,linkcolor=blue,citecolor=red]{hyperref}
\usepackage{caption}
\usepackage{subcaption}
\usepackage{url}
\usepackage{epstopdf}
\usepackage{pgf,tikz}
\usepackage{bbm}
\usepackage{extarrows}
\usetikzlibrary{arrows}
\usepackage{algorithm}
\usepackage{algorithmicx}
\usepackage{graphicx}       
\usepackage{amsmath}        
\usepackage{caption}
\usepackage{subcaption}
\usepackage{dsfont}
\usepackage{listings}
\usepackage{bbm}
\usepackage{color}
\usepackage{dsfont}
\usepackage{enumerate}
\usepackage{relsize}
\usepackage{xfrac}
\renewcommand{\Delta}{\triangle}
\newcommand{\bN}{\mathbb{N}}
\newcommand{\bE}{\mathbb{E}}

\definecolor{darkblue}{rgb}{0,0,0.7}

\definecolor{darkgreen}{rgb}{0.01,0.75,0.24}

\def \Ee[#1]{\mathcal{E}^{\text{{#1}}}}

\def\pa[#1,#2]{\frac{\partial {#1}}{\partial {#2}} }

\def\idom[#1,#2,#3]{\int_{#1}\hspace{1pt} {#2} \hspace{1pt} \text{d}{#3}}
\def\res[#1,#2]{\left.{#1}\right|_{#2}}

\def\var[#1,#2]{\langle \delta \mathcal{E}^{\text{{#1}}}({#2}),v\rangle}
\def\vars[#1,#2,#3]{\langle \delta^2\mathcal{E}^{\text{{#1}}}({#2})v,{#3}\rangle}
\def\vard[#1,#2,#3,#4]{\langle \delta\mathcal{E}^{\text{{#1}}}({#2})-\delta\mathcal{E}^{\text{{#3}}}({#4}),v\rangle}

\def\E{\mathbb{E}}

\newcommand{\cO}{\mathcal{O}}

\newcommand{\PP}{\mathbb{P}}

\newcommand{\blue}{\textcolor{black}}

\newcommand{\bZ}{\mathbb{Z}}
\newcommand{\bbR}{\mathbb{R}}
\newcommand{\bbP}{\mathbb{P}}
\newcommand{\interval}[1]{[\mspace{-2.7mu}[ #1 ]\mspace{-2.7mu}]}
\newcommand{\be}{\begin{equation}}
\newcommand{\en}{\end{equation}}
\newcommand{\ben}{\begin{equation*}}
\newcommand{\enn}{\end{equation*}}
\newcommand{\bea}{\begin{aligned}}
\newcommand{\ena}{\end{aligned}}

\def\ba#1\ena{\begin{align}#1\end{align}}
\def\ban#1\enan{\begin{align*}#1\end{align*}}

\theoremstyle{plain}
\newtheorem{thm}{Theorem}[section]

\newtheorem{definition}[thm]{Definition}

\newtheorem{corollary}[thm]{Corollary}

\newtheorem{assumption}[thm]{Assumption}

\newtheorem{proposition}[thm]{Proposition}
\newtheorem{remark}[thm]{Remark}

\numberwithin{equation}{section}

\begin{document}
%
%
\title[Improved efficiency of MLMC for SPDE through strong coupling]
{Improved efficiency of multilevel Monte Carlo \\
for stochastic PDE through strong pairwise coupling}
\author[N. K. Chada] {Neil K. Chada}
\address{Applied Mathematics and
Computational Science Program,
King Abdullah University of Science and Technology,
Thuwal 23955-6900, Kingdom of Saudi Arabia}
\email{neilchada123@gmail.com}
\author[H. Hoel] {H{\aa}kon Hoel}
\address{Chair of Mathematics for Uncertainty Quantification,
RWTH Aachen University, Aachen, Germany}
\email{hoel@uq.rwth-aachen.de}
\author[A. Jasra] {Ajay Jasra}
\address{Applied Mathematics and Computational Science Program,
King Abdullah University of Science and Technology,
Thuwal 23955-6900, Kingdom of Saudi Arabia}
\email{ajay.jasra@kaust.edu.sa}
\author[G. E. Zouraris] {Georgios E. Zouraris}
\address{Department of Mathematics and Applied Mathematics,
GR-700 13 Voutes Campus, Heraklion, Crete, Greece}
\email{georgios.zouraris@uoc.gr}
\begin{abstract}
  Multilevel Monte Carlo (MLMC) has become an important
  metho\-do\-lo\-gy in applied mathematics for reducing the
  computational cost of weak approximations. For many problems, it is
  well-known that strong pairwise coupling of numerical solutions in
  the multilevel hierarchy is needed to obtain efficiency gains.  In
  this work, we show that strong pairwise coupling indeed is also
  important when MLMC is applied to stochastic partial differential
  equations (SPDE) of reaction-diffusion type, as it can improve the
  rate of convergence and thus improve tractability.
For the MLMC method with strong pairwise coupling that was developed
and studied numerically on filtering problems in
[{\it Chernov et al., Numer. Math., 147 (2021), 71-125}],
we prove that the rate of computational efficiency is higher
than for existing methods.
We also provide numerical comparisons with alternative coupling ideas
on linear and nonlinear SPDE to illustrate the importance of this feature.
\end{abstract}
%
%
\maketitle
%
\section{Introduction}\label{sec:intro}
%
The efficiency of numerical methods is a very important topic for
practitioners that has lately seen a surge of interest in the field of
uncertainty quantification (UQ) \cite{TJS14,RCS13,DX10}. UQ seeks to
combine statistical and probabilistic techniques with traditional
numerical schemes to improve the modeling and the accuracy of estimates. 
Examples of applications include climate modeling,
subsurface flow, medical imaging and deep learning
\cite{APS21,TGF15,MW06}. A particular focus has been given on the
class of nume\-ri\-cal methods known as Monte Carlo (MC) methods,
which are used to solve problems incorporating elements of randomness
or uncertainty \cite{LPS14,RS13,DX10}, i.e., in stochastic computations.
%
%
One me\-tho\-do\-lo\-gy which has exhibited improved efficiency and a high
level of applicability, is multilevel Monte Carlo (MLMC).

MLMC is a numerical technique aimed at reducing the computational
cost of the Monte Carlo method. 
%
%
The methodology was first introduced by Heinrich
\cite{SH11} and extended and popularized by various works on diffusion
processes by Giles \cite{MBG08,MBG15}. The methodology of MLMC can be
viewed as a variance-reduction technique. Since these works, MLMC has
been applied in numerous areas, including stochastic
filtering~\cite{CJY20,FMS20,HLT16,JKL17,JKL18}, Markov chain Monte Carlo
(MCMC) \cite{CST13,DKS19} and partial differential equations with
random input arising in UQ~\cite{ABS13,CST13,HPS13}.  MLMC is based
upon a given problem, such as estimating an expectation at some
terminal time, w.r.t.~the law of a diffusion process, that requires a
discretization. For instance, in the diffusion case, this can be a
time-discretization based on the Euler method.  One then decomposes an
expectation w.r.t.~a law associated to a very precise discretization
into a telescoping sum of differences of expectations associated to
laws of increasingly coarse discretizations. The objective is then to
sample from coupled probability distributions associated to
consecutive discretized laws and to apply Monte Carlo at each summand
of the telescoping sum to achieve a variance reduction, relative to
using Monte Carlo at the finest discretization. The amount of
discretization refers to the level, in the acronym MLMC.

Despite the substantial advancements made with MLMC, the number of
applications and research papers on applying the methodology to
stochastic partial differential equations (SPDE)~\cite{AJ16} is
relatively small. Such examples include finite-difference solvers
with applications in mathematical finance~\cite{GR12}, and finite
element methods for parabolic SPDE~\cite{BL12,BLS13}. 
  There are open questions on the efficiency and scope of MLMC for SPDE which we 
  use as motivation for this work: is it possible to 
  improve the efficiency of MLMC through strong pairwise 
  coupling of numerical solutions of SPDE, and can that widen the scope of
  MLMC on SPDE to problems in higher dimensions and with lower-regularity driving noise?  

Our objective in this manuscript is to present a complexity study of an
alternative way to apply MLMC for SPDE, which can demonstrate
computational gains. This approach is based on the exponential Euler
method~\cite{EL19,HO10,JK09,LT13} and strong pairwise coupling of
solution realizations on different levels. The strong pairwise
coupling approach was introduced and studied experimentally in the
work on finite-dimensional Langevin SDE by M\"uller et
al.~\cite{MSS15} and extended to filtering methods for
(infinite-dimensional) SPDE by Chernov et al.~\cite{CHL20}. The
coupling idea is based on the exponential Euler
integrator~\cite{EL19,HO10,PGN11,LT13} for time-discretization of
reaction-diffusion type SDE/SPDE. For the finite-dimensional SDE in
\cite{MSS15}, strong coupling is shown to produce constant-factor efficiency gains 
in numerical experiments, whereas for the herein considered class of SPDE,
we show that strong coupling reduces the asymptotic rate
of growth in the computational cost. This indicates that 
strong coupling for MLMC can lead to more substantial
asymptotic efficiency gains for infinite-dimensional problems than for
finite-dimensional ones.

The main contribution of this work is to demonstrate the improvements
of the discussed coupling approach, for numerically solving
SPDE. This is presented in the standard-format cost-versus-error
result for exponential Euler MLMC in Theorem~\ref{thm:mainExpEuler}.  Specifically, our
findings suggest that in order to achieve $\mathcal{O}(\epsilon^{2})$ mean squared 
error (MSE) in a standard setting, we have to pay
$\mathcal{O}(\epsilon^{-2})$ in computational cost. This is a
reduction in cost compared to other existing methods, such as the 
Milstein MLMC method, for which the cost
is $\mathcal{O}(\epsilon^{-3})$, cf.~Theorem~\ref{thm:milsteinMLMC} and~\cite{BLS13}, 
and it is, to the best of our knowledge, the first theoretical result on the performance
of the exponential Euler MLMC for any nonlinear stochastic PDE. We also 
verify these gains numerically on two SPDE, one with a linear 
reaction term and one with a nonlinear one. 

The outline of this paper is as follows. In Section \ref{sec:MLMC} we
describe our model problem, which is a semilinear SPDE, and review
fundamental properties of the MLMC method. Section~\ref{sec:coup}
describes our proposed coupling method for MLMC and two alternative
methods. We also summarize the theoretical properties of our main MLMC
method in Theorem~\ref{thm:mainExpEuler}.  Numerical experiments on
various SPDE are conducted in Section~\ref{sec:num} to demonstrate the
improvement with the proposed coupling.  Finally, we conclude our
findings, and provide future areas of research, in
Section~\ref{sec:conc}. Required model assumptions are provided in the
Appendix.
\section{Background material}\label{sec:MLMC}
In this section we present and review the MLMC method applied to numerical discretizations of SPDEs. 
We first introduce the SPDE under consideration, and then review the approximation methods: a spectral 
Galerkin spatial discretization combined with either exponential Euler or Milstein discretization 
in time.
\subsection{Notation}
Let $T>0$ and let $(\Omega, \mathcal{F}, \mathbb{P})$ be a complete
probability space equipped with a filtration $\mathcal{F}_{t\in[0,T]}$.    
%
%
$H$ denotes a non-empty separable Hilbert space with inner product
$\langle \cdot, \cdot \rangle$, norm
$\|\cdot\|_{ H} = \sqrt{\langle \cdot, \cdot \rangle}$
and orthogonal basis $(e_n)_{n=1}^{\infty}$. 
$L^2(\Omega,H)$ denotes the associated Bochner--Hilbert space,
consisting of the set of strongly measurable maps $f:\Omega\to H$
such that
\[
\|f \|_{ L^2(\Omega,H)}^2  := \int_{\Omega} \|f(\omega)\|_{ H}^2 \, \mathrm{d}\bbP(\omega) < \infty
\]

Let $\bN := \{1,2,\ldots\}$, and for every $N\in\bN$, 
we introduce the finite-di\-men\-sion\-al
subspace $H^N:= \textrm{span} \{ e_n\mid n=1,\dots,N\}\subset H$
and the associated orthogonal projection operator 
$P_{ N} v := \sum^{ N}_{n=1} \langle v,e_n \rangle_{}\, e_n$
for $v\in H$.
For a (normally implicitly given) set $B$ and mappings $f,g:B \to [0,\infty)$, the notation $f \lesssim  g$ implies there exists a $C>0$ such that
$f(x) \le  C g(x)$ for all $x\in B$, and the notation $f\eqsim g$
means that both $f \lesssim g$ and $g\lesssim f$ hold. For multivariate positive-valued 
functions $f(x,y)$ and $g(x,y)$ for which it holds for some $C>0$ 
that $f(x,y) \le C g(x,y)$ for all $(x,y) \in \text{Domain}(f) = \text{Domain}(g)$, we 
write $f \lesssim_{(x,y)} g$ if confusion is possible.
And, similarly as above, $f \eqsim_{(x,y)} g$ means that 
$f \lesssim_{(x,y)} g$ and $g \lesssim_{(x,y)} f$. 
For $m,n \in \bZ$ with $m\le n$, we introduce the integer interval $\interval{m,n}:= [m,n]\cap\bZ$, 
and for $x\in\bbR$ we define $\lceil x \rceil:=\min\{n \in\bZ\mid n\ge x\}$.
\subsection{Problem setup}
We consider a semilinear stochastic partial differential
equation \cite{DZ92} of the form
\begin{equation}\label{eq:SPDE}
\begin{split}
dU_t=&\,\left(AU_t + f(U_t)\right)\,dt+dW_t
\quad\text{for}\quad t \in [0,T],\\
U_0=&\,u_0,
\end{split}
\end{equation}
where $A:D(A) \rightarrow H $ is a linear operator, $u_0\in H$ is a
random-valued initial condition, $f:H \rightarrow H$ is a
reaction term that in general is nonlinear, and $W_t$ is a
$Q$-Wiener process, cf.~\eqref{eq:qwiener}. 
A number of further assumptions are imposed for the problem,
which we have deferred to Appendix~\ref{subsec:model_assumptions}. Suffice it to say here that we do
assume that the linear operator is negative-definite and spectrally
decomposable in the considered basis:
\begin{equation}\label{eq:aOperatorFirst}
A v = - \sum_{k=1}^{\infty} \lambda_k \langle e_k, v\rangle  e_k
\end{equation}
and that the $Q$-Wiener process takes the form
\[
W(t,x) = \sum^{\infty}_{n=1}\sqrt{q_n}e_n w^n_t,
\]
where $(w^n_t)_{n=1}^{\infty}$ is a sequence of independent
scalar-valued Wiener processes. We note that the eigenbasis of the
operator $A$, $(e_n)_{n=1}^{\infty}$, also appears in the representation 
of the $Q$-Wiener process. The strictly positive sequence
$(\lambda_n)_{n=1}^\infty$ and the non-negative sequence $(q_n
)_{n=1}^\infty$ are further described in Appendix~\ref{subsec:model_assumptions}.

The mild solution to equation~\eqref{eq:SPDE} is 
an $H$-valued predictable process $(U_t)_{t\in[0,T]}$
satisfying
\begin{equation}\label{Mild1}
{\mathbb P}\left( \omega \in \Omega \, \Big| \,  \,U_t = e^{A t}u_0 + \int^t_0 e^{A(t-s)}f(U_s)ds +
\int^t_0e^{A(t-s)}dW_s\, \quad \forall t \in [0,T] \right)=1. 
\end{equation}
The general form of~\eqref{eq:SPDE} encapsulates numerous SPDE in
practice. We will introduce and numerically study some of these in
Section~\ref{sec:num}.

\subsection{Numerical methods}
\label{ssec:num_meth}
Numerical approximations of SPDE have traditionally been computed through
the use of finite difference methods and finite element methods
(FEM) \cite{AJ16,LPS14,ZK17}. For the relevance of this work, we
will utilize and discuss an alternative class of Galerkin-based solvers.
To motivate such an alternative class, we review some of these techniques
below.
%
%
\subsubsection{Continuous-time spectral Galerkin methods}
\label{subsec:spectralGalerkin}
For $N\in{\mathbb N}$, consider the Galerkin problem of solving the SPDE~\eqref{eq:SPDE} on the subspace $H^N$:
\begin{equation}\label{eq:SDE}
\begin{split}
dU^{ N}_t=&\, \big (A_{ N}U^{ N}_t+f_{ N}(U^{ N}_t) \big) \, dt
+dW^{ N}_t,\\
U^{ N}_0=&\,u_0^{ N}:=P_{ N}u_0,\\ 
\end{split}
\end{equation}
where $A_{ N}:=P_{ N}A$, $f_{ N}(v):=P_{ N}(f(v))$
and
$W^{ N}_t:= P_{ N} W_t=\sum^{ N}_{n=1}\sqrt{q_n}\,e_n\, w^n_t.$
It is well-known \cite{KP92} that~\eqref{eq:SDE} has a unique mild solution
given by
\begin{equation*}\label{eq:mild_s}
U^{ N}_t = e^{A_Nt} u^{ N}_0 +
\int^t_0e^{A_N(t-s)}f_{ N}(U^{ N}_s)\;ds
+\int^t_0 e^{A_N(t-s)}\;dW^{ N}_s.
\end{equation*}

We next discuss two time-discretizations of spectral Galerkin methods.

%
\subsubsection{The exponential Euler method}
\label{subsec:exponentialEuler}

For a given $J\in{\mathbb N}$, let $\Delta t =\frac{T}{J}$
and let $(t_j)_{j=0}^{ J}$ be the nodes of a uniformly spaced mesh 
of $[0,T]$, so that $t_j=j\,\Delta t$ for $j\in\interval{0,J}$.
Then, for given $N\in{\mathbb N}$, the exponential Euler
approximations $(V^{ N,J}_j)_{j=0}^{ J}\subset H^N$
of $(U^{ N}_{t_j})_{j=0}^J$ are defined by
\begin{equation}\label{eq:euler}
\begin{split}
V^{ N,J}_0:=&\,u_0^{ N},\\
V^{ N,J}_{j+1}=&\,e^{A_N\Delta t} V^{ N,J}_j
+A_{ N}^{-1}(e^{A_N\Delta t}-I)f_{ N}(V^{ N,J}_j)\\
&\quad+\int^{t_{j+1}}_{t_j}e^{A_N(t_{j+1}-s)}\,dW^{ N}_s
\qquad\forall j\in\interval{0,J-1},
\end{split}
\end{equation}
where  $A_{ N}^{-1}:H^N\to H^N$ denotes the
inverse operator of $A_{ N}$. Defining for $n \in \interval{1,N}$ the components of
$V^{ N,J}_j$ and $f_{ N}(\cdot)$ by
$V^{ N,J}_{j,n}:=\langle V^{ N,J}_j,e_n \rangle$,~\\
\mbox{$f_{{ N},n}(\cdot):=\langle f(\cdot),e_n\rangle$,} respectively, 
and recalling the spectral decomposition
of the operator $A$, we arrive at the recursive relation
\begin{equation}\label{eq:exp_euler}
V^{ N,J}_{j+1,n} = e^{-\lambda_n\Delta t}\,V_{j,n}^{ N,J}
+\frac{1-e^{-\lambda_n \Delta t}}{\lambda_n}
\,f_{{ N},n}(V^{ N,J}_j) + R_{j,n},
\end{equation}
where%
\begin{equation*}
R_{j,n} :=  \sqrt{q_n} \int_{t_j}^{t_{j+1}}
e^{-\lambda_n (t_{j+1} -s)}  dw_t^{n} \stackrel{\tt d}{=}
\mathcal{N}\left(0,\frac{q_n\,(1 -e^{-2\lambda_n\Delta t}) }{2\lambda_n}\right).
\end{equation*}
We recall from~\eqref{eq:aOperatorFirst} that $-\lambda_n$ denotes the $n$-th eigenvalue of the operator $A$,
see also Assumption~\ref{assum:oper} in Appendix~\ref{subsec:model_assumptions} for further details.
Convergence properties of the exponential Euler scheme has
been studied in~\cite{JK09}, where they demonstrate
strong convergence and highlight an improvement in the order of convergence
in time against traditional numerical schemes: 
\begin{proposition}[Jentzen and Kloeden~\cite{JK09}]\label{prop:exponentialEuler}
Let all assumptions in Appendix~\ref{subsec:model_assumptions} hold
for some $\phi \in (0,1)$ relating to the regularity of the Q-Wiener process. Then 
\begin{equation}\label{eq:KJ09}
\max_{j\in  \interval{0,J}}
\mathbb{E}\left[\|U_{t_j}-V^{ N,J}_j\|^2_{ H}\right]
\lesssim_{(N,J)}
\lambda^{-2\phi}_{ N} + \left(\frac{\log_2(J)}{J} \right)^{2},
\end{equation}
where $U$ is the mild solution \eqref{Mild1} of \eqref{eq:SPDE},
and $(V^{ N,J}_j)_{j=0}^{ J}$ denotes the exponential
Euler approximation of the mild solution, cf.~\eqref{eq:euler}.
\end{proposition}
We note that the first term on the RHS of \eqref{eq:KJ09} is related to the 
discretization in space and the second term 
is related to the discretization in time.

The performance of a numerical method will be measured by the  
computational cost required to reach a mean squared error (MSE) $\cO(\epsilon^2)$.
Computational cost refers to the number of computational operations, 
where we count each addition, subtraction, multiplication, division, and 
each draw of a Gaussian random variable as one computational operation. 
It follows from this definition that if $f =0$ in the SPDE~\eqref{eq:SPDE},
then no evaluation of the reaction term is needed and the computational 
cost of computing the final-time solution $V_J^{N,J}$ is $\cO(JN)$, 
as each time iteration of~\eqref{eq:euler} consists of $\cO(N)$
computational operations. When $f\neq 0$, however, the cost becomes a more complicated 
expression in general, and we make the following assumption
to simplify matters: 
\begin{assumption}[Cost of evaluating $f_N$]\label{assum:costF}
For any $N \in \bN$ and $V^N \in H^N$, the cost of
of evaluating $f_N(V^N)$ is $\cO(N\log_2(N))$. 
\end{assumption} 
When Assumption~\ref{assum:costF} holds, each evaluation of $f_N(V^{J,N}_j)$ 
costs $\cO(N \log_2(N))$, and this accumulates 
to 
\begin{equation}\label{eq:computationalCostScheme}
\mathrm{Cost}(V_{J}^{N,J}) = \cO\big(J \, \big(N +\mathrm{Cost}(f_N) \big)\, \big) = \cO(J N \log_2(N))
\end{equation}
for the final-time solution. 

\begin{remark} \label{rem:fft-spde}
  In the numerical scheme~\eqref{eq:euler}
  used in Proposition~\ref{prop:exponentialEuler} and in Assumption~\ref{assum:costF}
  it is tacitly assumed that the nonlinear reaction term $f_{N}(V^{N,J}_j)$ can be 
  evaluated exactly for any $N \in \bN$ and $V^{N,J} \in H^N$. 
  For nonlinear reaction terms $f$, this may however not be possible 
  in practice. In computations, we will employ the fast Fourier transform (FFT) to 
  approximate $f_N(V^{N,J}_j)$ on a uniform mesh with $N$ degrees of freedom in space
  for each iteration of~\eqref{eq:euler}, where we refer to~\cite[Section 6.3.1]{CHL20} 
  and~\cite{LPS14} for further details on this procedure. 
  The approximation of $f_N$ by FFT may introduce so-called aliasing errors 
  in the numerical solution, cf.~\cite[page 334]{jentzenRockner2015}.
  Aliasing errors are not covered in the mathematical analysis of this paper, 
  but we will include the cost of using FFT in the computational cost of 
  all numerical methods studied.
\end{remark}

\begin{remark}
Disregarding aliasing errors, Assumption~\ref{assum:costF} holds 
when computing $f_N(V^N)$ by FFT for Nemytskii-operators 
$f(U)(x) = g(U(x))$ where the mapping $g: \bbR \to \bbR$
additionally satisfies that one evaluation costs $\cO(1)$.
Then
\[
  f_N(V^{N}) = \mathrm{FFT}\big( g(V^N(x=0)), g(V^N(x=1/N)), \ldots, g(V^N(x=(N-1)/N)) \big)   
\] 
where the right-hand side costs $\cO(N \log_2(N))$ to evaluate. 
\end{remark}

%
%
%
%
%
\subsubsection{The Milstein method}
 The Milstein method has been extended from SDE to different 
 forms of parabolic SPDE with multiplicative noise in~\cite{barth2012milstein,jentzenRockner2015}. We will here consider the version developed 
 in~\cite{jentzenRockner2015}, since its scheme is easy to express in our 
 problem setting, and it is also easy to extend to an MLMC method. 
 Using the previously introduced discretization parameters in space and time and recalling that the operators $A$ and $Q$ share the same eigenspace, 
 the Milstein scheme~\cite[equation (28)]{jentzenRockner2015} 
 takes the form
 \[
  \begin{split}
    V_0^{N,J} &:= u_0^N\\
  V^{ N,J}_{j+1} 
  &= e^{A_N \Delta t} \left( V^{N,J}_j   + \Delta t\,f_N(V^{ N,J}_j) + W^{N}(t_{j+1})-W^{N}(t_{j}) \right)
  \end{split}
\]
for $j \in \interval{0,J-1}$. On the component level, the scheme is given by 
\begin{equation}\label{eq:Milstein}
V^{ N,J}_{j+1,n}= 
e^{-\lambda_n \Delta t} \left( 
  V_{j,n}^{N,J} + \Delta t f_{{ N},n}(V^{  N,J}_j) + \sqrt{q_n}\big(w^{n}(t_{j+1})-w^{n}(t_{j})\big) \right) 
\end{equation}
for $n \in \interval{1,N}$ and $j \in \interval{0,J-1}$. 

We next present strong convergence rates for the Milstein scheme restricted 
to the additive-noise setting. For extensions to various multiplicative-noise
settings, see~\cite{jentzenRockner2015,barth2012milstein}.  
\begin{proposition}[Jentzen and R\"ockner~\cite{jentzenRockner2015}]\label{prop:Milstein}
Let Assumption~\ref{assum:milstein} in Appendix~\ref{subsec:model-assumptions-milstein} be fulfilled for some values of $\phi \in (1/2,1)$, $\kappa \in [0,\phi)$ and $\theta \in [\max(\kappa,\phi-1/2) , \phi)$, where $\phi$ is the noise parameter introduced in Assumption~\ref{assum:noise}. Then it holds that
\begin{equation*}
\mathbb{E}\left[\left\|U_{T}-V^{ N,J}_{J}\right\|^2_{H}\right]
\lesssim_{(N,J)} \lambda^{-2\theta}_{N} + J^{-2\min(2(\theta - \kappa), \theta)},
\end{equation*}
where $U_t$ denotes the mild solution to the SPDE \eqref{eq:SPDE}
and  $(V^{ N,J}_j)_{j=0}^{ J}$ denotes the 
Milstein approximation to the mild solution, cf.~\eqref{eq:Milstein}.
\end{proposition}

\begin{proof}[Connecting the result to the literature]
Remark~\ref{rem:notationMilstein} in Appendix~\ref{subsec:model-assumptions-milstein}~associates our parameters $(\phi,\kappa,\theta)$ 
with corresponding ones in~\cite[Assumptions 1-4]{jentzenRockner2015}.
In our additive-noise setting with the operators $A$ and $Q$ having 
the same eigenspace, Proposition~\ref{prop:Milstein} follows from~\cite[Theorem 1]{jentzenRockner2015}.
\end{proof}

Even when disregarding the differences in the regularity assumptions, a comparison of the convergence rates for the exponential Euler and Milstein method is not straightforward since the rates for exponential Euler only depend on the single parameter $\phi$, while the rates of Milstein depend on two additional parameters, $\kappa$ and $\theta$. To simplify the comparison we impose additional constraints on the relationship between the 
parameters $\phi$, $\kappa$ and $\theta$:  
\begin{corollary}\label{cor:Milstein}
  For some value of $\phi \in (1/2,1)$, let Assumption~\ref{assum:milstein} in 
  Appendix~\ref{subsec:model-assumptions-milstein} be fulfilled for some $\kappa \in [0,\phi/2)$ 
  and all $\theta \in [\max(\kappa,\phi-1/2) , \phi)$. 
Then for any sufficiently small fixed $\delta>0$, it holds that  
\begin{equation*}
    \mathbb{E}\left[\left\|U_{T}-V^{ N,J}_{J}\right\|^2_{H}\right]
    \lesssim_{(N,J)} \lambda^{-2\phi +\delta}_{N} + J^{-2\phi + \delta}.
\end{equation*}
\end{corollary} 
\begin{proof}
  For any sufficiently small $\delta>0$, Assumption~\ref{assum:milstein} holds for some  $\kappa < \phi/2 -\delta/4$ and $\theta_\delta :=\phi -\delta/2$. 
  Noting that
  \[
    \min(2(\theta_\delta - \kappa), \theta_\delta) = \theta_\delta = \phi - \delta/2,
  \] 
  the result follows from Proposition~\ref{prop:Milstein}.
\end{proof} 

For a fixed value of $\phi \in (1/2,1)$, the additional constraints imposed on $\kappa$ and $\theta$ in Corollary~\ref{cor:Milstein} are likely to present the Milstein method in a good light, as they produce the highest possible convergence rates attainable from Proposition~\ref{prop:Milstein}.
Comparing the convergence of exponential Euler in Proposition~\ref{prop:exponentialEuler}
with Milstein in Corollary~\ref{cor:Milstein}, the methods have essentially the same rate in space, but exponential Euler has a higher rate in time. Note further that the rates only apply to Milstein when $\phi >1/2$, while they apply to exponential Euler method for any $\phi \in (0,1)$. But one should also keep in mind that the Milstein method applies to a wider range of reaction terms $f$ than exponential Euler, since Assumption~\ref{assum:milstein}
is more relaxed than Assumption~\ref{assum:func}. When comparable, the lower convergence rate for Milstein leads to a poorer performance for the Milstein MLMC method than the exponential Euler MLMC method in low-regularity settings, when $\phi <3/4$, cf.~Theorems~\ref{thm:mainExpEuler} 
and~\ref{thm:milsteinMLMC}. See also Section~\ref{sec:num} for numerical evidence that exponential Euler outperforms Milstein when $\phi \approx 1/2$.

\subsection{The multilevel Monte Carlo method}
%
The expectation of an $H$-valued random variable $U$ 
is often approximated by the standard Monte Carlo estimator 
\begin{equation*}
E_{ M}\left[U\right]:=
\frac{1}{M}\,\sum^{ M}_{m=1}U^{(m)},
\end{equation*}
where the samples $U^{(1)}, U^{(2)},\ldots, U^{(M)} \sim \PP_U$ are
independently drawn random variables and $E_M[U]$ consequently denotes the sample average 
estimator using $M$ i.i.d.~draws of $U$. When it is computationally
costly to draw samples of $U$, variance-reduction techniques may improve the
efficiency through reducing the statistical error of the estimator.
The multilevel Monte Carlo (MLMC) method is an extension of 
standard Monte Carlo that draws pairwisely coupled random variables $\{(U^{\ell-1,C}, U^{\ell,F})\}_{\ell=0}^L$, where $U^{\ell-1,C}$ denotes the coarse random variable on resolution level $\ell$, and $U^{\ell,F}$ the fine random variable on level $\ell$. Pairwise coupling of $(U^{\ell-1,C}, U^{\ell,F})(\omega)$ means that $U^{\ell-1,C}(\omega)$ and $U^{\ell,F}(\omega)$ are generated using the same driving noise $W_t(\omega)$ (to be elaborated on in the next section). We further impose that
\begin{equation}\label{eq:telescoping_sum}
U^{-1,C} := 0 \in H \qquad \text{and} \qquad \E{U^{\ell,C}} = \E{U^{\ell,F}} \quad \forall \ell \in \bN_0,
\end{equation}
so that the weak approximation on resolution level $L \in \bN$
can be represented as a telescoping sum of expectations:  
\begin{equation}\label{eq:telescoping2}
\E{U} \approx \E{U^{L,F}} \stackrel{\eqref{eq:telescoping_sum}}{=} \sum_{\ell=0}^{L} \E{ U^{\ell,F}- U^{\ell-1,C}}.
\end{equation}
By approximating each of the $L+1$ expectations in the telescoping sum by a sample average, we obtain the MLMC estimator:
\begin{equation}\label{eq:mlmcEstimator}
\begin{split}
E_{ \textrm{M\!L}}\left[ U\right]:=&\,\sum_{\ell=0}^{ L} E_{ M_\ell}
\left[U^{\ell,F}- U^{\ell-1,C}\right]\\
=&\,\sum_{m=1}^{ M_0}\frac{U^{0,F,(m)}}{M_{0}}
+\sum^{ L}_{\blue{\ell=1}} \sum^{ M_{\ell}}_{m=1}
\frac{U^{\ell,F,(m)} - U^{\ell-1,C,(m)}}{M_{\ell}}.
\end{split}
\end{equation}
Here, $(U^{\ell-1,C,(m)},  U^{\ell,F,(m)})$ denotes the $\mathbb{P}_{(U^{\ell-1,C},U^{\ell,F})}$-distributed $m$-th sample on level $\ell$, and all samples on all resolution levels are independent, meaning that all random variables in the sequence $\{(U^{\ell-1,C,(m)}, U^{\ell,F,(m)})\}_{\ell,m}$ are independent. A near-optimal calibration of the parameters $L \in \bN$ and $(M_{\ell})_{\ell=0}^{ L}\subset\bN$ is obtained through minimizing the mean squared error for a given computational cost, cf.~\cite{MBG08} and 
Theorem~\ref{thm:VMLMC}. 
The MLMC estimator achieves variance reduction over standard Monte Carlo when the coupled random variables $U^{\ell-1,C}$ and $U^{\ell,F}$
are sufficiently correlated, cf.~Condition $(ii)$ in Theorem~\ref{thm:VMLMC} below.

\par 
One way to assess the performance of Monte Carlo methods is through the MSE. 
The following theorem describes the cost versus
error of the MLMC methodology for $H$-valued random variables:
\begin{thm}\label{thm:VMLMC}
Assume that the telescoping-sum properties~\eqref{eq:telescoping_sum} 
hold and that there exists positive constants 
$\alpha,\beta,\gamma$
such that $\alpha \geq\frac{\min(\beta,\gamma)}{2}$ and
\begin{itemize}
\item[(i)] $\left\|\bE\left[U^{\ell,F}-U\right]\right\|_{ H}
  \lesssim\,2^{-\alpha\,\ell}$,

  \item[(ii)] $V_\ell:=\mathbb{E}\left[\left\|U^{\ell,F}-U^{\ell-1,C}
      \right\|^2_{ H}\right]\lesssim\,2^{-\beta\,\ell}$,

  \item[(iii)]
    $C_\ell := \mathrm{Cost}(U^{\ell-1,C}, U^{\ell,F}) \lesssim\, 2^{\gamma\,\ell}$.
    
\end{itemize}
Then for any $\epsilon \in(0,1)$ and $L := \lceil \log_2(1/\epsilon)/\alpha \rceil$,
there exists a sequence $(M_{\ell})_{\ell=0}^{ L}\subset\bN$ such
that
\begin{equation*}
\mathrm{MSE}={\mathbb E}\left[\big\|E_{ \mathrm{M\!L}}\left[
U\right]-{\mathbb E}\left[U\right]\big\|^2_{ H}
\right]\lesssim\,\epsilon^2,
\end{equation*}
and
\begin{equation}\label{eq:mlmcCost}
\mathrm{Cost(MLMC)}:= \sum_{\ell=0}^L M_{\ell}\, C_{\ell}
\lesssim
\begin{cases} \epsilon^{-2}, \quad &\mathrm{if} \ \beta
>\gamma,\\ \epsilon^{-2}( \log \epsilon)^2, \quad
&\mathrm{if} \ \beta =\gamma,\\
\epsilon^{-2-\frac{(\gamma-\beta)}{\alpha}}, \quad
&\mathrm{if} \ \beta <\gamma.
\end{cases}
\end{equation}
\end{thm}
The proof of this result is a straightforward extension of
the original theorem presented by Giles \cite{MBG08} for weak
approximations of stochastic differential equations.
\begin{proof}
Let
\begin{equation}\label{eq:ml1}
  M_{\ell}:=\left \lceil \epsilon^{-2} \sqrt{\frac{V_\ell}{C_\ell}}
    \sum_{j=0}^L\sqrt{V_j C_j} \right\rceil \qquad \ell \in \interval{0,L},
\end{equation} 
where $V_0 := \mathbb{E}\left[\|U^{0}\|_{ H}^2\right]$.
By the telescoping-sum property
\[
  \bE[E_{ \mathrm{M\!L}}[U]]= \sum_{\ell=0}^L \E{U^{\ell,F} - U^{\ell-1,C}}
  \stackrel{\eqref{eq:telescoping_sum}}{=} \E{U^{L,F}},
\]
the representation~\eqref{eq:mlmcEstimator} and the independence 
of the samples \\ $\{(U^{\ell-1,C,(m)}, U^{\ell,F,(m)})\}_{\ell,m}$, we 
obtain that
\begin{equation*}
\begin{split}
\mathbb{E}\Big[\big\|E_{ \mathrm{M\!L}}[U]& -\mathbb{E}[U]\big\|^2_{ H}\Big]\\
&= \|\E{U} - \E{U^{ L}}\|_{H}^2 + \E{\left\|E_{ \mathrm{M\!L}}[U] -{\mathbb E}[U^{ L}]\right\|^2_{ H}}\\
&\lesssim   2^{-2\alpha L} + \E{ \left\| 
  \sum_{m=1}^{ M_0}\frac{U^{0,F,(m)} - \bE[U^0]}{M_{0}} \right\|^2 }\\
   &\quad +\E{  \left\|\sum^{L}_{\ell=1}\sum^{ M_{\ell}}_{m=1}
 \frac{U^{\ell,F,(m)} - U^{\ell-1,C,(m)} - \bE[U^{\ell,F} - U^{\ell-1,C}]}{M_{\ell}}  \right\|^2}\\
  &\lesssim\,\sum_{\ell=0}^{ L} \frac{V_\ell}{M_\ell} +2^{-2\alpha L} \lesssim \,\epsilon^2.
\end{split} 
\end{equation*}
By assumptions $(ii)$ and $(iii)$, we obtain that
\[
  \begin{split}
  \sum_{\ell=0}^LC_{\ell}\,  M_{\ell}\  &\stackrel{\eqref{eq:ml1}}{\le}
  \sum_{\ell=0}^LC_{\ell} \left(\epsilon^{-2} \sqrt{\frac{V_\ell}{C_\ell}}
  \sum_{j=0}^L\sqrt{V_j C_j} +1 \right) \\
  &\lesssim  \epsilon^{-2}  \left(\sum_{j=0}^L \sqrt{V_j C_j}\right)^2 + \underbrace{C_L}_{\lesssim 2^{\gamma L} } \\
& \lesssim \epsilon^{-2}  \left(\sum_{j=0}^L \sqrt{V_j C_j}\right)^2 + \epsilon^{- \gamma/\alpha}\\
& \lesssim  \begin{cases} 
  \epsilon^{-2}  & \text{if} \quad \beta > \gamma \\
  L^2 \epsilon^{-2} + \epsilon^{-2}& \text{if} \quad \beta = \gamma\\
  \epsilon^{-2-\frac{(\gamma-\beta)}{\alpha}} + \epsilon^{-\gamma/\alpha}& \text{if} \quad 
  \beta< \gamma.
\end{cases}
\end{split}
\]
For the last inequality, 
the assumption $\alpha \ge \min(\beta, \gamma)/2$ implies that 
that $\gamma/\alpha \le 2$ when $\beta \ge \gamma$ and $\beta/\alpha \le 2$
when $\beta\le \gamma$ (so that $2 + (\gamma-\beta)/\alpha \ge \gamma/\alpha$),
and inequality~\eqref{eq:mlmcCost} follows.
\end{proof}
\begin{remark}
The theorem \blue{also} applies in settings where one replaces
$V_\ell$ in Theorem \ref{thm:VMLMC} (ii) by
${\widetilde V}_\ell:={\mathbb E}\left[
\left\| U^{\ell,F}-U^{\ell-1,C}
-\bE\left[U^{\ell,F}-U^{\ell-1,C}\right]
\right\|^2_{ H}\right]$,
and for some problems this may improve the rate $\beta>0$.
Practically, however, there may be little to gain by replacing $V_\ell$ 
by $\widetilde V_\ell$ as weak approximations of $U^{\ell,F} - U^{\ell-1,C}$
can be much more intractable than strong approximations, cf.~\cite{LP18}.
\end{remark}
\section{Multilevel Monte Carlo methods for SPDE}\label{sec:coup}
In this section we describe two MLMC methods that are based on extending 
the two numerical schemes in Section~\ref{ssec:num_meth} to the MLMC setting. 
To better illustrate the importance of strong coupling and the loss of
accuracy due to damping, we also propose a third MLMC method which is an extension of a 
modified form of the exponential Euler method that only is 
exponential in the drift-term.  
We will employ the following notation for the multilevel hierarchy of
discretized solutions: On level $\ell \ge 0$, let 
$N_{\ell}\eqsim N_0\,2^{\nu \ell}$ for given $N_0\in\bN$ and $\nu >0$
denote a sequence of spatial resolutions, and let 
$J_\ell:=J_0\,2^\ell$ for a given $J_0 \in \bN$ denote a sequence of time resolutions. 
In a notation that suppresses details on the pairwise coupling, we let
$U^{\ell,F}_j := V^{ N_\ell, J_\ell}_j$ 
denote the fine numerical solution of a given spectral Galerkin method on 
level $\ell$ at time $t_j^{\ell} :=j\,\Delta t_\ell$ for $j\in\interval{0,J_\ell}$,
computed on the subspace $H^{N_\ell}$ using the time step 
$\Delta t_{\ell}:=\frac{T}{J_\ell}$.
And $U^{\ell-1,C}_j := V^{ N_{\ell-1}, J_{\ell-1}}_j$ denotes the coupled coarse numerical solution 
on level $\ell$ at time $t_j^{\ell-1} :=j\,\Delta t_{\ell-1}$ for $j\in\interval{0,J_{\ell-1}}$ computed on the subspace $H^{N_{\ell-1}}$ with time step $\Delta t_{\ell-1}:=\frac{T}{J_{\ell-1}}$.

To discuss the quality of a pairwise coupling, let us first introduce some terminology. When a coupling satisfies  
\[
  \E{U^{\ell,F}_{j}} = \E{U^{\ell, C}_{j}} \quad \forall j \in \interval{0,J_{\ell}}
\]
for all $\ell\ge0$, we say that the coupling is \emph{weakly correct},
and when it additionally satisfies  
\[
  U^{\ell,F}_{j}(\omega) = U^{\ell,C}_{j}(\omega)  \quad \forall (\omega, j) \in \Omega \times \interval{0,J_\ell}
\]
for all $\ell\ge0$, we say that it is a \emph{pathwise correct coupling}. 
From the construction of the multilevel estimator in Section~\ref{sec:MLMC}, we see that weakly correct coupling is needed to obtain the crucial telescoping sum in the MLMC estimator,~cf.~\eqref{eq:telescoping_sum} and~\eqref{eq:telescoping2}, and that weakly correct coupling thus ensures consistency for the MLMC estimator. Pathwise correct coupling is on the other hand not necessary to obtain consistency, and there are many examples of performant MLMC methods that only are weakly correct, cf.~\cite{giles2014antithetic,hoel2020mlenkf}. Pathwise correct coupling 
is however often an easy way to ensure the needed weakly correct coupling.  

To achieve high performance, the pairwise coupling must be weakly correct and produce a high convergence rate $\beta$ for the strong error, cf.~Theorem~\ref{thm:VMLMC}. We will refer to a coupling that achieves a high rate $\beta$ in comparison to alternative approaches as a \emph{strong coupling}. To be more precise for the particular SPDE considered in this work, we introduce the notion of strong diffusion coupling: 
\begin{definition}[Strong diffusion coupling (SDC)]
  Consider a weakly correct coupling sequence of spectral-Galerkin numerical solutions of the \\ $\{(U^{\ell-1,C}, U^{\ell,F})\}_{\ell\ge0}$ of the SPDE~\eqref{eq:SPDE} with no reaction term, $f=0$ (the stochastic heat equation). Recall further that a coupled pair of solutions is defined on time meshes of different resolutions:  
  \[
    U^{\ell-1,C}_j = U^{\ell-1,C}(j \Delta t_{\ell-1}) \in H^{N_{\ell-1}} \quad \text{for} \quad j \in \interval{0, J_{\ell-1}}
  \] 
  and 
  \[
    U^{\ell,F}_j = U^{\ell,F}(j \Delta t_{\ell}) \in H^{N_{\ell}} \quad \text{for} \quad j \in \interval{0, J_{\ell}},
  \] 
  with $\Delta t_{\ell-1} = 2 \Delta t_\ell$. We say that the coupling is a strong diffusion coupling if it holds for all $\ell \ge 0$ that  
  \begin{equation*}\label{eq:strongDiffusionCoupling}
   P_{N_{\ell-1} } U^{\ell,F}_{2j} = U^{\ell-1,C}_{j} \quad \forall j \in \interval{0, J_{\ell-1}}.
  \end{equation*}
\end{definition}  

For the stochastic heat equation, an SDC is thus an exact coupling of $U^{\ell-1,C}$ to $U^{\ell,F}$ on the subspace $H^{N_{\ell-1}}$. This is of course the strongest possible coupling one can achieve (for the given problem), and we will see later that the exponential Euler MLMC method indeed is the only among the three we consider whose coupling is SDC. Although our theory and numerical experiments both indicate a connection between SDC and strong couplings more generally 
when $f$ is non-zero-valued, it is not clear how far this extends. To best of our knowledge, it is an open problem to describe coupling strategies for $H$-valued stochastic processes that are weakly correct and maximize the convergence rate of the strong error $\beta$.

We next extend the exponential Euler method and the Milstein method to the MLMC setting.


%
\subsection{Exponential Euler MLMC method}\label{ssec:exponentialEulerMLMC}

\blue{This MLMC method was first introduced and analyzed for the linear
reaction-term setting in~\cite[Section 5.4.1]{CHL20}.
Since then the method has been applied to the SPDE \eqref{eq:SPDE} with linear 
reaction term for problems arising in Bayesian computation. These include stochastic filtering \cite{JLX21} and Markov chain Monte Carlo \cite{JKL18}, with an extension to multi-index Monte Carlo.}


We consider the pairwisely coupled solutions $(U^{\ell-1,C},U^{\ell,F})$ 
that both are solved by the numerical scheme~\eqref{eq:exp_euler}
with the respective initial conditions
\[
U^{\ell-1,C}_0 = P_{N_{\ell-1}}u_0 \quad \text{and} \quad 
U^{\ell,F}_0 = P_{N_\ell}u_0.   
\]

For the fine solution, the $n$-th component 
of two iterations of the scheme~\eqref{eq:exp_euler} 
at time $t^{\ell}_{2j} = 2j \Delta t_{\ell}$ 
takes the form 
\begin{equation}\label{eq:fineSolExponential1}
U_{2j+1,n}^{\ell,F} = e^{-\lambda_n \Delta t_{\ell}} U_{2j,n}^{\ell,F}
+ \frac{1- e^{- \lambda_n \Delta t_{\ell} }}{\lambda_n}
f_{{ N_{\ell}},n}(U_{2j}^{\ell,F} ) + R_{2j,n}^{\ell,F} ,
\end{equation}
and
\begin{equation}\label{eq:fineSolExponential2}
U_{2j+2,n}^{\ell,F} = e^{- \lambda_n \Delta t_{\ell} } U_{2j+1,n}^{\ell,F}
+\frac{1- e^{-\lambda_n \Delta t_{\ell} }}{ \lambda_n}
f_{{ N_{\ell}},n}(U_{2j+1}^{\ell,F} ) + R_{2j+1,n}^{\ell,F},
\end{equation}
for $(j,n) \in \interval{0, J_{\ell-1}-1} \times \interval{1,N_{\ell}}$
and with
\begin{equation}\label{eq:couplFine}
R_{k,n}^{\ell,F} =\sqrt{q_n}\,
\int_{t^{\ell}_{k}}^{t^{\ell}_{k+1}}
e^{-\lambda_n (t^{\ell}_{k+1}-s)}\;dw^{n}_s
\stackrel{{\tt d}}{=} 
N\left(0, q_n \frac{1- e^{- 2\lambda_n \Delta t_{\ell} }}{2 \lambda_n}\right)
\end{equation}
for $(k,n) \in \interval{0,J_{\ell}-1} \times \interval{1,N_\ell}$. 

The coupled coarse solution uses the time step $\Delta t_{\ell-1} = 2\,\Delta t_{\ell}$,
and one iteration at time $t_{j}^{\ell-1} = j\Delta t_{\ell-1} = 2j\Delta t_{\ell} = t_{2j}^{\ell}$
takes the form  
\begin{equation}\label{eq:coarseSolExponential}
U_{j+1,n}^{\ell-1,C}= e^{- \lambda_n \Delta t_{\ell-1} } U_{j,n}^{\ell-1,C}
+\frac{1-e^{-\lambda_n\,\Delta t_{\ell-1}}}{ \lambda_n}
f_{{ N_{\ell-1}},n}(U^{\ell-1,C}_{j,n})+ R_{j,n}^{\ell-1,C},
\end{equation}
where
\[
R_{j,n}^{\ell-1,C} =   \sqrt{q_n} \int_{t_j^{\ell-1}}^{t_{j+1}^{\ell-1} } e^{-\lambda_n (t_{j+1}^{\ell-1} -s)} dw^{n}_s.
\]
The pairwise coupling $U^{\ell-1,C}_{j+1,n} \leftrightarrow U^{\ell,F}_{2j+2,n}$ is obtained
through coupling the driving noise $R^{\ell-1,C}_{j,n}
\leftrightarrow (R_{2j,n}^{\ell,F},R_{2j+1,n}^{\ell,F})$. 
By~\eqref{eq:couplFine}, we have that
\begin{equation*}
\begin{split}
R_{j,n}^{\ell-1,C}=&\,\sqrt{q_n}\,
\int_{t_j^{\ell-1}}^{t_{j+1}^{\ell-1} } e^{-\lambda_n(t^{\ell-1}_{j+1} - s)} dw^{n}_s\\
=&\,e^{-\lambda_n \Delta t_{\ell}}\,\sqrt{q_n}
\,\int_{t^{\ell}_{2j}}^{t^{\ell}_{2j +1} } 
e^{-\lambda_n (t^{\ell}_{2j+1} -s)} dw^{n}_s
+\sqrt{q_n}\,\int_{t^{\ell}_{2j+1}}^{t^{\ell}_{2j+2}} 
e^{-\lambda_n (t^{\ell}_{2j+2} -s)} dw^{n}_s,\\
\end{split}
\end{equation*}
which yields
\begin{equation}\label{eq:correct}
R_{j,n}^{\ell-1,C} = e^{-\lambda_n \Delta t_{\ell}} R_{2j,n}^{\ell,F}
+ R_{2j+1,n}^{\ell,F}
\qquad \forall\,(j,n)
\in \interval{0, J_{\ell}-1} \times \interval{1, N_{\ell}}.
\end{equation}

To summarize, given the coupling $U^{\ell-1,C}_{j,n} \leftrightarrow U^{\ell,F}_{2j,n}$ at some time $t_j^{\ell-1}$,
the coupling at the next time is obtained by 
generating the fine-solution noise $(R_{2j,n}^{\ell,F},R_{2j+1,n}^{\ell,F})$ 
and coupling it to the coarse-solution noise by formula~\eqref{eq:correct}. The next-time solution $U^{\ell-1,C}_{j+1,n}$ is computed by~\eqref{eq:coarseSolExponential}
with $R_{j,n}^{\ell-1,C}$ as input, and
$U^{\ell,F}_{2j+2,n}$ is computed 
by~\eqref{eq:fineSolExponential1} and~\eqref{eq:fineSolExponential2}
with $(R_{2j,n}^{\ell,F},R_{2j+1,n}^{\ell,F})$ as input.

\begin{remark}\label{rem:coupling}
We note from the above that 
\[
  R^{\ell-1,C}_{j,n}(\omega) =  \sqrt{q_n} \int_{t_j^{\ell-1}}^{t_{j+1}^{\ell-1} } e^{-\lambda_n (t_{j+1}^{\ell-1} -s)} dw^{n}_s(\omega) = R^{\ell-1,F}_{j,n}(\omega), 
\]  
and since $U^{\ell-1,C}$ and $U^{\ell-1,F}$ are solved using the same numerical scheme,
the coupling is pathwise correct.
Let us further note that if $f=0$, then 
the linearity of the problem and~\eqref{eq:correct} imply
that the coupling is an SDC: 
\begin{equation}\label{eq:pathwise-true-f0}
  U^{\ell,F}_{2j, n}(\omega)  = U^{\ell-1,C}_{j, n}(\omega)  \quad
\forall (j,n) \in \interval{0, J_{\ell-1}} \times \interval{1,N_{\ell-1}}.
\end{equation}
This can be verified by induction: assume~\eqref{eq:pathwise-true-f0}
holds for some $j\in \interval{0,J_{\ell-1}-1}$ (it holds for $j=0$ by definition). And using the numerical schemes for the respective methods with $f=0$,
we obtain that
\[
  \begin{split}
  U^{\ell,F}_{2j+2, n} &\stackrel{\eqref{eq:fineSolExponential2}}{=} 
  e^{- \lambda_n \Delta t_{\ell} } U_{2j+1,n}^{\ell,F}+ R_{2j+1,n}^{\ell,F} \\
  &\stackrel{\eqref{eq:fineSolExponential1}}{=}
 e^{- \lambda_n 2\Delta t_{\ell} } U_{2j,n}^{\ell,F}
 + e^{- \lambda_n \Delta t_{\ell} } R_{2j,n}^{\ell,F} + R_{2j+1,n}^{\ell,F}\\
 &\stackrel{\eqref{eq:correct}}{=}  e^{- \lambda_n \Delta t_{\ell-1} } U_{j,n}^{\ell-1,C}
+ R_{j,n}^{\ell-1,C}\\
&\stackrel{ \eqref{eq:coarseSolExponential}}{=} U_{j+1,n}^{\ell-1,C}.
\end{split} 
\]
Since the exponential Euler MLMC method is SDC, we expect it to perform very efficiently when $f=0$, and Theorem~\ref{thm:mainExpEuler} shows that the coupling is strong also for more general reaction terms. 
\end{remark}

For showcasing the importance of strong pairwise coupling, and as a transition between exponential Euler MLMC and Milstein MLMC, we next consider a 
slightly altered form of the exponential Euler method with explicit integration of the 
It\^o integral. 
\subsection{Drift-exponential Euler MLMC method}
We consider the drift-exponential Euler scheme 
\begin{equation*}\label{eq:drift-exp-euler}
  \begin{split}
  V^{ N,J}_0:=&\,u_0^{ N},\\
  V^{ N,J}_{j+1}=&\,  e^{A_N\Delta t} V^{ N,J}_j + A_{ N}^{-1}(e^{A_N\Delta t}-I)f_{ N}(V^{ N,J}_j)\\
  & + e^{A_N\Delta t} \big( W^N(t_{j+1}) - W^N(t_j) \big) \Big)
  \qquad\forall j\in\interval{0,J-1}.
\end{split}
\end{equation*}
This is a mix of exponential Euler and Milstein, as the approximation of the drift terms agree with the exponential Euler scheme and the
approximation of the It\^o integral agrees with the Milstein scheme. 

When extending this scheme to an MLMC method, a similar argument as in Section~\ref{ssec:exponentialEulerMLMC} yields that 
two iterations of the fine solution in the pairwise couple $(U^{\ell-1,C}, U^{\ell,F})$ 
takes the form 
\begin{equation}\label{eq:fineSolDriftExponential1}
  U_{2j+1,n}^{\ell,F} = e^{-\lambda_n \Delta t_{\ell}} U_{2j,n}^{\ell,F}
  + \frac{1- e^{- \lambda_n \Delta t_{\ell} }}{\lambda_n}
  f_{{ N_{\ell}},n}(U_{2j}^{\ell,F} ) + \widetilde R_{2j,n}^{\ell,F} ,
  \end{equation}
  and
  \begin{equation}\label{eq:fineSolDriftExponential2}
  U_{2j+2,n}^{\ell,F} = e^{- \lambda_n \Delta t_{\ell} } U_{2j+1,n}^{\ell,F}
  +\frac{1- e^{-\lambda_n \Delta t_{\ell} }}{ \lambda_n}
  f_{{ N_{\ell}},n}(U_{2j+1}^{\ell,F} ) + \widetilde R_{2j+1,n}^{\ell,F},
  \end{equation}
  for $(j,n) \in \interval{0, J_{\ell-1}-1} \times \interval{1,N_{\ell}}$
  and with
  \begin{equation}\label{eq:couplDriftExponentialFine}
  \widetilde R_{k,n}^{\ell,F} = \sqrt{q_n} \, e^{-\lambda_n \Delta t_{\ell}} \Big(w^{n}(t_{k+1}^\ell)- w^{n}(t_k^\ell)\Big)
  \end{equation}
  for $(k,n) \in \interval{0,J_{\ell}-1} \times \interval{1,N_\ell}$. 
  
  The coupled coarse solution takes the form  
  \begin{equation}\label{eq:coarseSolDriftExponential}
  U_{j+1,n}^{\ell-1,C}= e^{- \lambda_n \Delta t_{\ell-1} } U_{j,n}^{\ell-1,C}
  +\frac{1-e^{-\lambda_n\,\Delta t_{\ell-1}}}{ \lambda_n}
  f_{{ N_{\ell-1}},n}(U^{\ell-1,C}_{j,n})+ \widetilde R_{j,n}^{\ell-1,C},
  \end{equation}
  where
  \[
  \widetilde R_{j,n}^{\ell-1,C} =   \sqrt{q_n} e^{-\lambda_n \Delta t_{\ell-1}} \big(w^{n}(t_{j+1}^{\ell-1}) - w^{n}(t_j^{\ell-1}) \big).
  \]
  Recalling that $t_{j}^{\ell-1} = j \Delta t_{\ell-1} = 2j \Delta t_\ell = t^{\ell}_{2j}$, 
  we obtain the pairwise coupling of $U^{\ell-1,C}_{j+1,n} \leftrightarrow U^{\ell,F}_{2j+2,n}$
  through coupling the driving noise: 
  \begin{equation*}\label{eq:driftExplicitMLMC}
    \begin{split}
  \widetilde R_{j,n}^{\ell-1,C} &= \sqrt{q_n} e^{-2\lambda_n \Delta t_{\ell}} \big(w^{n}(t_{2j+2}^{\ell}) - w^{n}(t_{2j}^{\ell}) \big)= e^{-\lambda_n \Delta t_{\ell}}\big( \widetilde R_{2j,n}^{\ell,F} + \widetilde R_{2j+1,n}^{\ell,F}\big).
    \end{split} 
  \end{equation*}
Since $\widetilde R_{j,n}^{\ell-1,C}(\omega) = \widetilde R_{j,n}^{\ell-1,F}(\omega)$ and
$U^{\ell-1,C}$ and $U^{\ell-1,F}$ are solved using the same numerical method, it follows that 
that the drift-exponential Euler MLMC method also is pathwisely correctly coupled.
However, it is not SDC, since when $f=0$ we obtain by~\eqref{eq:fineSolDriftExponential1} and~\eqref{eq:fineSolDriftExponential2} that for $n \in \interval{1,N_{\ell-1}}$,  
\[
  \begin{split}
U^{\ell,F}_{2,n} &= e^{- \lambda_n 2\Delta t_{\ell} } U^{\ell,F}_{0,n}
+  \widetilde R_{1,n}^{\ell,F} + e^{- \lambda_n \Delta t_{\ell}} \widetilde R_{0,n}^{\ell,F} \\
&= \underbrace{e^{- \lambda_n \Delta t_{\ell-1} } U^{\ell,C}_{0,n} +  \widetilde R_{0,n}^{\ell-1,C}}_{= U^{\ell-1,C}_{2,n}} 
+ (1 - e^{- \lambda_n \Delta t_{\ell}} )\widetilde R_{1,n}^{\ell,F} \neq U^{\ell-1,C}_{2,n}.
  \end{split}
\]
The term $ (1 - e^{- \lambda_n \Delta t_{\ell}} )\widetilde R_{1,n}^{\ell,F}$ 
is an error in the coupling that is introduced by explicit integration of the It\^o integral. This leads to an artificial smoothing of the numerical solution, as is illustrated
by the numerical examples in Section~\ref{sec:num}. 

\subsection{Milstein MLMC method}
We consider the pairwise coupling of the coarse and fine Milstein 
solutions on level $\ell$ with respective initial conditions 
\[
  U^{\ell-1,C}_0 = P_{N_{\ell-1}}u_0 \quad \text{and} \quad U^{\ell,F}_0 = P_{N_\ell}u_0.
\]
Two iterations of the fine solution 
takes the form
\begin{equation*}\label{eq:fineSolMilstein1}
  \begin{split}
  U^{\ell,F}_{2j+1,n} &= 
  e^{-\lambda_n \Delta t_\ell} \left( 
    U^{\ell,F}_{2j,n}  + \Delta t_\ell f_{N_{\ell},n}(U^{\ell,F}_{2j,n}) + \sqrt{q_n}\big(w^{n}(t_{2j+1}^{\ell})-w^{n}(t_{2j}^{\ell}) \big) \right) \\
    &= e^{-\lambda_n \Delta t_\ell} \left( 
      U^{\ell,F}_{2j,n}  + \Delta t_\ell f_{N_{\ell},n}(U^{\ell,F}_{2j,n}) \right) +  \widehat R_{2j,n}^{\ell,F}
\end{split}  
\end{equation*}
and 
\begin{equation*}\label{eq:fineSolMilstein2}
  \begin{split}
  U^{\ell,F}_{2j+2,n} &= 
  e^{-\lambda_n \Delta t_\ell} \left( 
    U^{\ell,F}_{2j+1,n}  + \Delta t_{\ell} f_{N_{\ell},n}(U^{\ell,F}_{2j+1,n}) + \sqrt{q_n}\big(w^{n}(t_{2j+2}^{\ell})-w^{n}(t_{2j+1}^{\ell}) \big) \right) \\
    &= e^{-\lambda_n \Delta t_\ell} \left( 
      U^{\ell,F}_{2j+1,n}  + \Delta t_{\ell} f_{N_{\ell},n}(U^{\ell,F}_{2j+1,n}) \right) +  \widehat R_{2j+1,n}^{\ell,F}
  \end{split}
\end{equation*}
for $(j,n) \in \interval{0,J_{\ell-1} -1} \times \interval{1,N_\ell}$ with 
\[
  \widehat R_{k,n}^{\ell,F} := \sqrt{q_n} e^{-\lambda_n \Delta t_\ell} \left(w^{n}(t_{k+1}^{\ell})-w^{n}(t_{k}^{\ell})\right), \; \quad 
  (k,n) \in \interval{0,J_\ell-1} \times \interval{1, N_\ell}.
\] 

One iteration of the coarse solution takes the form 
\begin{equation*}\label{eq:coarseSolMilstein}
  \begin{split}
  U^{\ell-1,C}_{j+1,n} &= 
  e^{-\lambda_n \Delta t_{\ell-1}} \left( 
    U^{\ell-1,C}_{j,n}  + \Delta t_{\ell-1} f_{N_{\ell-1},n}(U^{\ell-1,C}_{j,n}) \right)\\
    &    + \underbrace{\sqrt{q_n} e^{-\lambda_n \Delta t_{\ell-1}} \left(w^{n}(t_{j+1}^{\ell-1})-w^{n}(t_{j}^{\ell-1}) \right)}_{=: \widehat R^{\ell-1,C}_{j,n}}
  \end{split}
\end{equation*}
for $(j,n) \in \interval{0,J_{\ell-1} -1} \times \interval{1,N_{\ell-1}}$.

And we obtain the same coupling as for the drift-exponential Euler method:
\begin{equation*}\label{eq:couplingMilstein}
\widehat R^{\ell-1,C}_{j,n}  = \sqrt{q_n} e^{-\lambda_n 2\Delta t_{\ell}} \left(w^{n}(t_{2j+2}^{\ell})-w^{n}(t_{2j}^{\ell}) \right) 
=  e^{-\lambda_n \Delta t_{\ell}}\left(\widehat R^{\ell,F}_{2j,n} + \widehat R^{\ell,F}_{2j+1,n} \right). 
\end{equation*}
By a similar argument as for the previous MLMC method, this is a pathwise correct coupling, 
but it is not SDC. 

In summary, we have presented three different MLMC methods where only the coupling for the exponential Euler MLMC method is SDC. This particularly means that the exponential Euler MLMC method outperforms the other methods when $f=0$, and later comparisons of the strong convergence rate
$\beta$ for the two methods in Theorems~\ref{thm:mainExpEuler} and~\ref{thm:milsteinMLMC} and in the numerical experiments show that the outperformance is broader.  

\subsection{MLMC for SPDE}

In this section, we present cost versus error results for the exponential Euler- and Milstein
MLMC methods.

We recall that the computational cost of one simulation of a numerical method
is defined by the computational effort needed, cf.~\eqref{eq:computationalCostScheme},
and that under Assumption~\ref{assum:costF}, it holds for all three spectral 
Galerkin methods we consider that   
\[
\text{Cost}(U^{\ell,F}(T)) = \text{Cost}(U^{\ell,F}_{J_\ell}) = \text{Cosst}(V^{N_\ell,J_\ell}_{J_\ell}) \eqsim  J_\ell N_\ell \log_2(N_\ell),
\] 
and
\begin{equation}\label{eq:costBound}
\text{Cost}( U^{\ell-1,C}(T),  U^{\ell,F}(T) )
\eqsim  \text{Cost}( U^{\ell-1,C}_{J_{\ell-1}}) + \text{Cost}(U^{\ell,F}_{J_\ell} ) \eqsim  J_\ell N_\ell \log_2(N_\ell).  
\end{equation}

\blue{
We will consider weak approximations of Banach-space-valued 
quantities of interest (QoI) of the following form:
\begin{definition}[Admissible QoI] \label{def:admissibleQoI}
  Let $K$ be a Banach space equipped with the norm $\|\cdot \|_K$ and let $\varphi: H \to K$ 
  be a strongly measurable and uniformly Lipschitz continuous QoI. We say that 
  such a QoI is admissible if the cost of evaluating the mapping 
  satisfies that 
  \[
  \sup_{v \in H^N} \text{Cost}(\varphi(v)) \lesssim N.  
  \]
  \end{definition}
}

We are ready to state the main result of this work.

\begin{thm}[Exponential Euler MLMC]
  \label{thm:mainExpEuler}
  Consider the SPDE~\eqref{eq:SPDE} for a linear operator $A$ with
  $\lambda_n \eqsim n^2$ \blue{and let $\varphi:H \rightarrow K$ 
  be an admissible QoI, in the sense of Definition~\ref{def:admissibleQoI}.} If all assumptions in
  Appendix~\ref{subsec:model_assumptions} hold for some $\phi \in (0,1)$ and Assumption~\ref{assum:costF} holds, then the pathwise correctly coupled exponential Euler MLMC method with
  \begin{equation}\label{eq:resolutionExpEuler}
    J_\ell = 2^\ell J_0 \quad \text{ and } \quad N_\ell \eqsim N_0 2^{\ell/(2\phi)} 
  \end{equation}
  satisfies
  \begin{itemize}
  \item[(i)] $\big\| \bE \big[ \blue{\varphi}(U^{\ell,F}(T,\cdot))  -  \blue{\varphi}(U(T, \cdot)) \big]\big\|_{\blue{K}} \lesssim
    (\ell+1) 2^{-\ell}$.
  \item[(ii)] $V_\ell :=  \mathbb{E}\Big[\| \blue{\varphi}(U^{\ell,F}(T,\cdot))  - \blue{\varphi}( U^{\ell-1,C}(T, \cdot))\|^2_{\blue{K}} \Big] \lesssim (\ell +1)^2 \; 2^{-2\ell}$.
  \item[(iii)] $C_\ell := \mathrm{Cost}\big(\blue{\varphi}(U^{\ell-1,C}(T)), \blue{\varphi}(U^{\ell,F}(T)) \big) \blue{\lesssim} (\ell +1) 2^{(1+1/(2\phi))\ell}$.
  \end{itemize}
And for any sufficiently small $\epsilon>0$
and $L := \lceil \log_2\big( \log_2(1/\epsilon)/\epsilon \big) \rceil$,
there exists a sequence $\{M_{\ell}(\epsilon)\}_{\ell=0}^L \subset \bN$ such
that
\begin{equation}\label{eq:mlmc_mse}
\mathrm{MSE} = \mathbb{E}\big[\big\| E_{ \mathrm{M\!L}}
[\blue{\varphi}(U(T,\cdot))]-\mathbb{E}[ \blue{\varphi}(U(T, \cdot)) ] \big\|^2_{\blue{K}} \big] \lesssim  \epsilon^{2},
\end{equation}
and 
\begin{equation}\label{eq:mlmcCost_correc}
  \begin{split}
\mathrm{Cost(MLMC)}  &:= \sum_{\ell=0}^L M_\ell C_\ell\\ 
  &\lesssim  \begin{cases} 
  \epsilon^{-2}  & \text{if} \quad \phi \in  (1/2,1) \\
  \epsilon^{-2} \big(\log_2(1/\epsilon)\big)^{5}  & \text{if} \quad \phi = 1/2\\
  \epsilon^{-1 - 1/(2\phi)} \big(\log_2(1/\epsilon)\big)^{2+1/(2\phi)} & \text{if} \quad  \phi \in (0,1/2).
\end{cases}
\end{split}
\end{equation}
\end{thm}
\begin{proof}
  
  \blue{
 Let us show that the $K$-valued random variables $\varphi(U^{\ell,F}(T,\cdot))$ and $\varphi(U^{\ell,C}(T,\cdot))$ are well-defined. The Lipschitz continuity of the mapping $\varphi$
implies that
\[
  \begin{split}
  \|\varphi(U^{\ell,F}(T,\cdot))\|_K &\le  \|\varphi(U^{\ell,F}(T,\cdot)) - \varphi( 0)\|_K 
  + \|\varphi( 0)\|_K\\
  & \le  C_{\varphi} \|U^{\ell,F}(T,\cdot)\|_H + \|\varphi(0)\|_{K}, 
  \end{split}
  \] 
  where $C_{\varphi}>0$ denotes the Lipschitz constant for $\varphi$.
  It follows that $\varphi(U^{\ell,F}(T,\cdot)) \in L^2(\Omega, K)$ for all $\ell \ge 0$,
  and we similarly also have that $\varphi(U^{\ell,C}(T,\cdot)) \in L^2(\Omega, K)$. 
}

\blue{We note that the numerical resolution sequences are set according to~\eqref{eq:resolutionExpEuler} 
to balance the error from space- and time-discretization in Proposition~\ref{prop:exponentialEuler}.}
A pair of correctly coupled solutions $U^{\ell,F, (m)}$
and $U^{\ell-1,C,(m)}$ can be viewed as exponential Euler solutions using the
same driving $Q-$Wiener process $W^{(m)}$ on levels $\ell$ and $\ell-1$, respectively. Consequently,
\begin{equation*}
\begin{split}
  \mathbb{E}&\left[\|\blue{\varphi}(U^{\ell,F}(T,\cdot))-\blue{\varphi}(U^{\ell-1,C}(T,\cdot))\|^2_{\blue{K}}\right]\\
& \leq \blue{C_{\varphi}^2} \mathbb{E}\left[\|U^{\ell,F}(T,\cdot)-U^{\ell-1,C}(T,\cdot)\|^2_{\textcolor{black}{{H}}}\right]\\
&\leq\,2 \blue{C_{\varphi}^2}\mathbb{E}\left[ \| U^{\ell,F}(T,\cdot)-U(T,\cdot)\|_{ H}^2 + \|U(T,\cdot) - U^{\ell-1,C}(T, \cdot)\|^2_{ H} \right]\\
&= 2 \blue{C_\varphi^2} \mathbb{E}\left[ \| V^{N_\ell,J_\ell}_{J_\ell}-U(T,\cdot)\|_{ H}^2 + \|U(T,\cdot) - V^{N_{\ell-1}, J_{\ell-1}}_{J_{\ell-1}}\|^2_{ H} \right]\\
& \hspace*{-0.1cm}\stackrel{\eqref{eq:KJ09} }{\lesssim}\, \lambda_{N_{\ell}}^{-2\phi} + \left(\frac{\log_2(J_\ell)}{J_\ell}\right)^2 \\
&\eqsim\, N_\ell^{-2\phi} + \Big(\ell+1\Big)^2 2^{-2\ell} \\
&\eqsim\, (\ell+1)^2 \, 2^{-2\ell},
\end{split}
\end{equation*}
\blue{where we used Lipschitz continuity for the first inequality. 
This verifies rate (ii). Since $\{\varphi(U^{\ell,F}(T,\cdot))\}_\ell$ is a Cauchy sequence in $L^2(\Omega,K)$ with limit $\varphi(U(T,\cdot))$, we have that 
\[
  \begin{split}
\E{\blue{\varphi}(U(T, \cdot)) -\blue{\varphi}(U^{\ell,F}(T,\cdot))}
&= \E{\sum_{j = \ell}^\infty \blue{\varphi}(U^{j+1, F}(T, \cdot)) -\blue{\varphi}(U^{j,F}(T,\cdot))}\\
&= \E{\sum_{j = \ell}^\infty \blue{\varphi}(U^{j+1, F}(T, \cdot)) -\blue{\varphi}(U^{j,C}(T,\cdot))}.
  \end{split}  
\]
In the last equality we used that the coupling is pathwise correct: $U^{j,F}(T,\cdot) = U^{j,C}(T,\cdot)$, cf.~Remark~\ref{rem:coupling}. Rate (i) follows from  
\[
  \begin{split}
\big\| \bE \big[ \blue{\varphi}(U(T, \cdot)) -\blue{\varphi}(U^{\ell,F}(T,\cdot))  \big]\big\|_{\blue{K}} & \le 
\sum_{j=\ell}^\infty \bE \big[  \|\blue{\varphi}(U^{j+1, F}(T, \cdot)) -\blue{\varphi}(U^{j,C}(T,\cdot))  \|_K \big] \\
&\stackrel{(ii)}{\lesssim} \sum_{j=\ell+1}^{\infty} (j+1) 2^{-j}\\
&\lesssim (\ell+1)2^{-\ell}.
\end{split} 
\]
}
\blue{Rate (iii) follows by~\eqref{eq:costBound} and Definition~\ref{def:admissibleQoI}.
Introducing the following number-of-samples-per-level sequence
\begin{equation}\label{eq:ml2}
M_\ell =  \left \lceil \epsilon^{-2} \sqrt{\frac{V_\ell}{C_\ell}}
\sum_{j=0}^{ L}\sqrt{V_j C_j} \right\rceil
\quad\ell\in\interval{0,L},
\end{equation}
and noting that $L = \lceil \log_2(\log_2(1/\epsilon)/\epsilon) \rceil \eqsim \log_2(1/\epsilon)$,
we obtain~\eqref{eq:mlmc_mse} by a similar argument as in the proof of Theorem~\ref{thm:VMLMC}:
\[
  \begin{split}
\mathbb{E}\big[\big\| E_{ \mathrm{M\!L}} [\blue{\varphi}(U(T,\cdot))] -\mathbb{E}[ \blue{\varphi}&(U(T, \cdot)) ] \big\|^2_{\blue{K}} \big] \\  
&= \mathbb{E}\Big[\left\|E_{ \mathrm{M\!L}}[ \blue{\varphi}(U(T,\cdot))] -{\mathbb E}[\varphi(U^{L,F}(T,\cdot))]\right\|^2_{\blue{K}}\Big]\\
& \quad + \|\mathbb{E}[\blue{\varphi}(U(T,\cdot))] 
- \blue{\varphi}(U^{L,F}(T,\cdot))]\|_{\blue{K}}^2\\
&\lesssim   \sum_{\ell=0}^{ L} \frac{V_\ell}{M_\ell} + 
L^2 2^{-2L}\\
&\lesssim \epsilon^{2} \frac{\max( \sum_{\ell=0}^L \sqrt{V_\ell C_\ell}, \, 1 )}{\max(\sum_{j=0}^L \sqrt{V_j C_j}, \, 1)} +L^2 \frac{\epsilon^2}{(\log(\epsilon))^2}   \\ 
&\eqsim \,\epsilon^2.
  \end{split}
  \] 
}
  For the computational cost, we have that 
  \[
  \begin{split}
  \sum_{\ell=0}^LC_{\ell}\,  M_{\ell}\ 
  &\lesssim  \epsilon^{-2}  \left(\sum_{j=0}^L \sqrt{V_j C_j}\right)^2 + C_L \\
& \lesssim \epsilon^{-2}  \left(\sum_{j=0}^L \sqrt{V_j C_j}\right)^2 + \left( \frac{\epsilon}{ \log_2(1/\epsilon)} \right)^{-( 1+ 1/(2\phi) )},
\end{split}
\]
and~\eqref{eq:mlmcCost_correc} follows from using
$V_jC_j \lesssim (j+1)^3\, 2^{\,j\, (1/(2\phi)-1)}$
when bounding the squared sum from above.
\end{proof}

\begin{remark}
A general framework for (MLMC) methods for reaction-diffusion type
SPDE in the setting of $\phi \ge 1/2$ and for numerical methods with
a strong convergence rate $1/2$ was first developed
in~\cite{BLS13}. When $\phi=1/2$ and $\gamma = 2$, the MSE
$\mathcal{O}(\epsilon^{2-\delta})$ was achieved at the computational
cost $\cO(\epsilon^{-3})$ for that method in~\cite[Theorem
4.4]{BLS13} compared to a cost $\mathcal{O}(\epsilon^{-2})$ for our 
exponential Euler MLMC method. This is however not a fair performance 
comparison, since~\cite{BLS13} 
was developed for more general SPDE with multiplicative noise 
and for which the operators $A$ and $Q$ need not share eigenbasis,
while our method is tailored to the additive-noise setting
with $A$ and $Q$ sharing eigenbasis, cf.~Appendix~\ref{subsec:model_assumptions}.
\end{remark}

We state a similar cost-versus-error result for the
MLMC Milstein method with pathwise correctly pairwise coupling. 
\begin{thm}[Milstein MLMC]
  \label{thm:milsteinMLMC}
  Consider the SPDE~\eqref{eq:SPDE} for a linear operator $A$ with
  $\lambda_n \eqsim n^2$, let the assumptions in Corollary~\ref{cor:Milstein} hold for 
some $\phi \in (1/2,1)$ and let Assumption~\ref{assum:costF} hold. 
\blue{Let $\varphi:H \rightarrow K$ 
  be an admissible QoI, in the sense of Definition~\ref{def:admissibleQoI}.} 
Then the pathwise correctly coupled Milstein MLMC method with 
\begin{equation}\label{eq:resolutionMilstein}
  J_\ell = 2^\ell J_0 \quad \text{ and } \quad N_\ell \eqsim N_0 2^{ \ell/2},
\end{equation}
satisfies for any fixed $\delta>0$ that
\begin{itemize}
\item[(i)] $\big\| \bE \big[ \blue{\varphi}(U^{\ell,F}(T,\cdot))   -  \blue{\varphi}(U(T, \cdot)) \big]\big\|_H \lesssim 2^{-(\phi-\delta/2) \ell}$.
\item[(ii)] $V_\ell :=  \mathbb{E}\Big[\| \blue{\varphi}(U^{\ell,F}(T,\cdot))  -  \blue{\varphi}(U^{\ell-1,C}(T, \cdot))\|^2_H \Big] \lesssim  2^{-(2\phi-\delta) \ell}$.
\item[(iii)] $C_\ell := \mathrm{Cost}\big( \blue{\varphi}(U^{\ell-1,C}(T)) ,  \blue{\varphi}(U^{\ell,F}(T)) \big) \eqsim (\ell +1) 2^{3\ell/2}$.
\end{itemize}

And for any sufficiently small fixed $\delta >0$ and any sufficiently small $\epsilon>0$, there exist an 
$L(\epsilon) \in \bN$ and a sequence $\{M_{\ell}(\epsilon)\}_{\ell=0}^L \subset \bN$ such
that 
\begin{equation}\label{eq:mlmc_mse_milstein}
  \mathrm{MSE} = \mathbb{E}\big[\big\| E_{ \mathrm{M\!L}}
  [ \blue{\varphi}(U(T,\cdot))]-\mathbb{E}[ \blue{\varphi}(U(T, \cdot)) ] \big\|^2_H \big] \lesssim  \epsilon^{2},
  \end{equation}
  at the cost  
  \begin{equation}\label{eq:mlmcCost_milstein}
    \begin{split}
  \mathrm{Cost(MLMC)} &:= \sum_{\ell=0}^L M_\ell C_\ell 
    \lesssim  \begin{cases} 
    \epsilon^{-2}  & \text{if} \quad \phi \in  (3/4,1) \\
    \epsilon^{-2(1+\delta)}   & \text{if} \quad \phi = 3/4\\
    \epsilon^{-3(1 +\delta)/(2\phi)}  & \text{if} \quad  \phi \in (1/2,3/4).
  \end{cases}
  \end{split}
  \end{equation}
\end{thm}

\begin{proof}
We set the numerical resolution sequences by~\eqref{eq:resolutionMilstein}
to balance the error from space- and time-discretization in Corollary~\ref{cor:Milstein},
and, \blue{since the Milstein MLMC method is pathwise correctly coupled}, 
the rates (i), (ii) and (iii) can be verified as in the proof of Theorem~\ref{thm:mainExpEuler}.
 
To prove the error and cost results, we relate the rates in (i), (ii) and (iii) to those in 
Theorem~\ref{thm:VMLMC}: for any $\delta>0$, it holds that 
\begin{equation}\label{eq:alphabetagamma}
  \alpha=\phi-\delta/2,  \quad \beta =2\phi -\delta  \quad \text{and} \quad \gamma = 3/2+\delta
\end{equation} 

For the case $\phi \in (3/4,1)$ it holds for sufficiently small $\delta>0$ that $\beta> \gamma$, and the results~\eqref{eq:mlmc_mse_milstein} 
and~\eqref{eq:mlmcCost_milstein} follow from Theorem~\ref{thm:VMLMC}. 

For the case $\phi \in (1/2,3/4]$, we again apply Theorem~\ref{thm:VMLMC}
to our rates $(\alpha,\beta,\gamma)$ in~\eqref{eq:alphabetagamma}
to conclude that~\eqref{eq:mlmc_mse_milstein} is fulfilled at the cost
\[
  \mathrm{Cost(MLMC)} \lesssim \epsilon^{-2-\frac{(\gamma-\beta)}{\alpha}} = 
  \epsilon^{-\frac{3/2 +\delta}{\phi -\delta/2}},
\]
and taking $\delta>0$ sufficiently small, it holds that 
\[
  \frac{3/2 +\delta}{\phi -\delta/2} \le \frac{3(1 + \delta)}{2\phi}.
\]
\end{proof}

Comparing Theorem~\ref{thm:mainExpEuler} with Theorem~\ref{thm:milsteinMLMC},
we expect exponential Euler MLMC to asymptotically outperform Milstein MLMC 
when the colored noise has low regularity, meaning when $\phi < 3/4$.
%

\section{Numerical examples}
\label{sec:num}

In this section, we numerically test the exponential
Euler MLMC method against the drift-exponential- and Milstein MLMC methods.
We study two reaction-diffusion SPDE, 
one with a linear reaction term and one with a trigonometric
one. To showcase the superior performance 
of exponential Euler in settings with low-regularity colored noise, 
we consider one setting with $\phi \approx 1/2$ (this is a low-regularity
setting for the Milstein method) and we numerically confirm the 
theoretical result that the exponential Euler MLMC and the Milstein MLMC 
perform similarly when $\phi \approx 3/4$, cf.~Theorems~\ref{thm:mainExpEuler} 
and~\ref{thm:milsteinMLMC}.

\par 
For our numerical experiments we consider the general form
of semilinear SPDE
\[
dU_t =\,\left(AU_t + f(U_t)\right)\;dt+dW_t, \quad t \in [0,T],
\]
with initial triangular-wave initial condition 
\[
u_0(x) = \begin{cases} 2x, & x \in\left[0,\frac{1}{2}\right], \\
  2(1-x), & x \in\left(\frac{1}{2},1\right].
\end{cases}
\]
Furthermore we specify our space $H = L^2(0,1)$ with Fourier basis functions $e_n(x) = \exp(i 2n\pi x)$ 
for $n \in \mathbb{Z}$, and the final time is set to $T=1/2$.
We consider a linear operator $A:D(A) \rightarrow H$ defined as
\begin{equation*}
A =-\sum_{n \in \mathbb{Z}} \lambda_n\,
\langle \cdot, e_n\rangle \,e_n,
\end{equation*}
with eigenvalues $(\lambda_n)_{n \in \mathbb{Z}}$, given as
\begin{equation*}
  \lambda_n = \begin{cases}
    1 & \mathrm{if} \quad n = 0,\\
    \frac{(2 n \pi)^2}{5}  & \mathrm{if} \quad  n \in \mathbb{Z}\setminus\{0\}.
    \end{cases}
\end{equation*}
We note that the triangular-wave initial condition satisfies the following regularity condition: $u_0 \in H_{3/4-\delta}$ for any $\delta>0$.

For $f:H \rightarrow H$, we consider the two different reaction terms which
are presented in Table \ref{table:mappings}. Both belong
to the class of Nemytskii operators, cf.~\cite{LPS14}.
\begin{table}[h!]
\begin{center}
    \begin{tabular}{| l | p{3cm} |}
    \hline
     $f(U)(x)$ &\textbf{Reaction term} \\ \hline \hline
      $U(x)$ & Linear  \\ \hline
     $2(\sin(2 \pi U(x) ) + \cos( 2\pi U(x)))$ & Trigonometric \\ 
    \hline
    \end{tabular} \bigskip \smallskip
        \caption{Different reaction terms $f(U):H \to H$ tested in the numerical experiments.}
    \label{table:mappings} 
\end{center}
\end{table}

The driving noise $dW$ is a $Q-$Wiener process~\eqref{eq:qwiener}
with 
\[
  q_n := \frac{1}{4}\lambda_n^{-2b} \qquad  \text{for}  \quad n\in\mathbb{Z},
\]
for two different values of $b$: the low-regularity setting $b=1/4$, and the smoother setting $b=1/2$. 
In connection with Assumption~\ref{assum:noise}, we note that 
\[
\sum_{n \in \mathbb{Z}} (\lambda_n)^{2\phi -1} q_n =   
\frac{1}{4} + \frac{2\pi^2}{5} \sum_{n =1}^\infty n^{4(\phi -b) -2} < \infty \iff \phi < 1/4 +b. 
\]
It consequently holds that that $\phi =(1/4+b)-\delta$ for any $\delta>0$, and for simplicity,
we will refer to the parameter values for $\phi$ as $\phi(b=1/4) = 1/2-$ and $\phi(b=1/2) =3/2-$, respectively. 
When Theorems~\ref{thm:mainExpEuler} and~\ref{thm:milsteinMLMC} apply, we expect exponential Euler MLMC to outperform Milstein MLMC when $\phi <3/4$, and that the methods perform similarly when $\phi > 3/4$. 

Note however that some of our numerical studies are purely experimental, as neither of the theorems apply to all problem settings we consider. Theorem~\ref{thm:mainExpEuler} only applies to the linear reaction term, 
because the trigonometric reaction term has no Fr{\'e}chet derivative that belongs to $L(H)$, and this violates Assumption~\ref{assum:func}.  
We do however believe the regularity assumptions in Proposition~\ref{prop:exponentialEuler} can be relaxed so that it also applies to the trigonometric reaction term, but, to the best of our knowledge, it is 
an open problem to prove this. 

For the Milstein method, on the other hand, Assumption~\ref{assum:milstein} does hold whenever $\phi > 1/2$ and $\kappa > 1/4$, with Fr{\'e}chet derivatives $f'(\cdot) = 4\pi (\cos(2\pi \cdot) -\sin(2\pi \cdot))$ and $f''(\cdot) = -8\pi^2f(\cdot)$. (This can be verified using the definition of Fr{\'e}chet derivatives and that $L^\infty(0,1) \subset H_{\kappa}$.) But Theorem~\ref{thm:milsteinMLMC} only applies when $\phi >1/2$.  

\subsection{Numerical estimates of the convergence rate $\beta$}

\begin{figure}[h!] \center
  \includegraphics[width=0.49\textwidth]{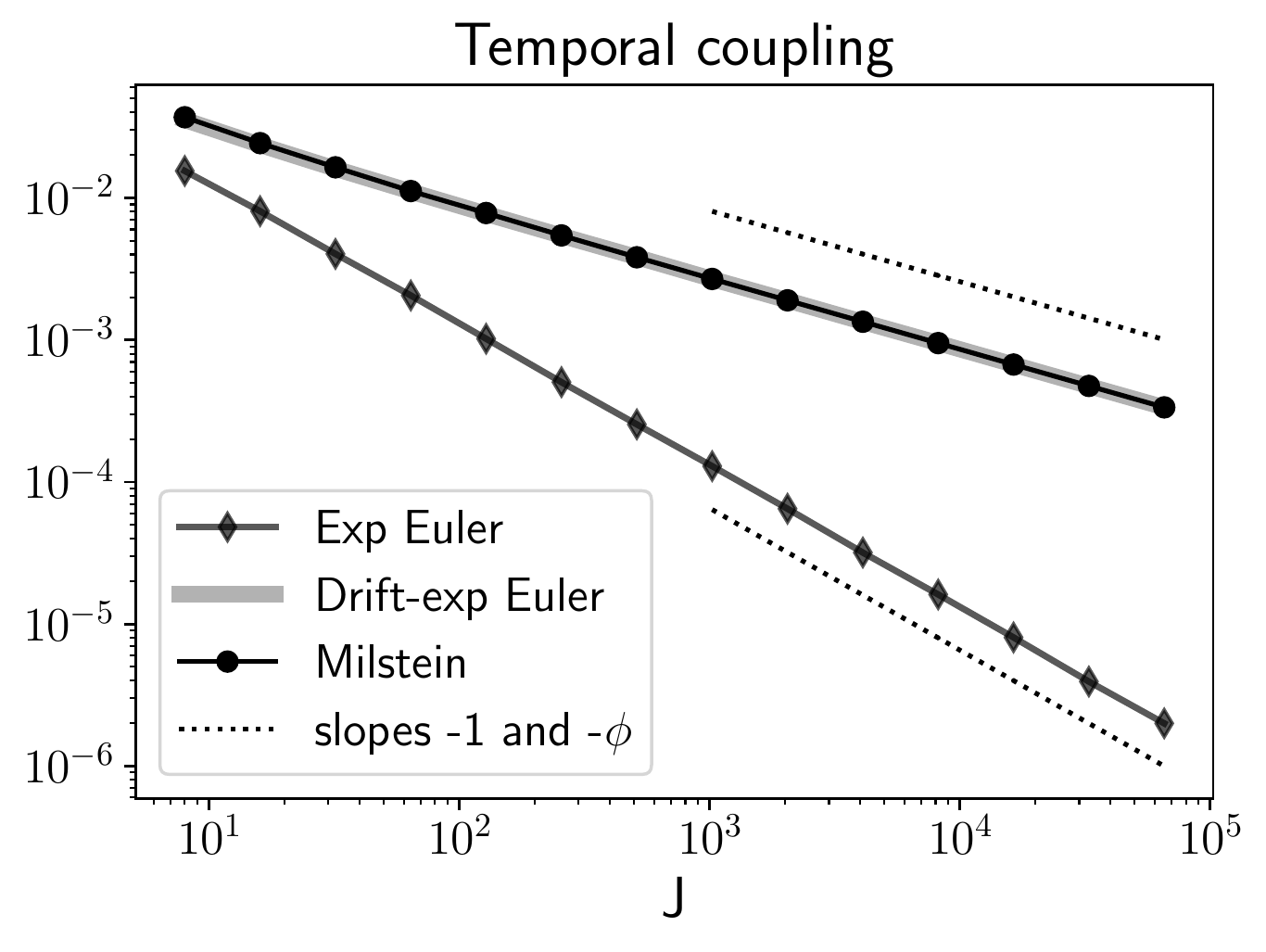}
  \includegraphics[width=0.49\textwidth]{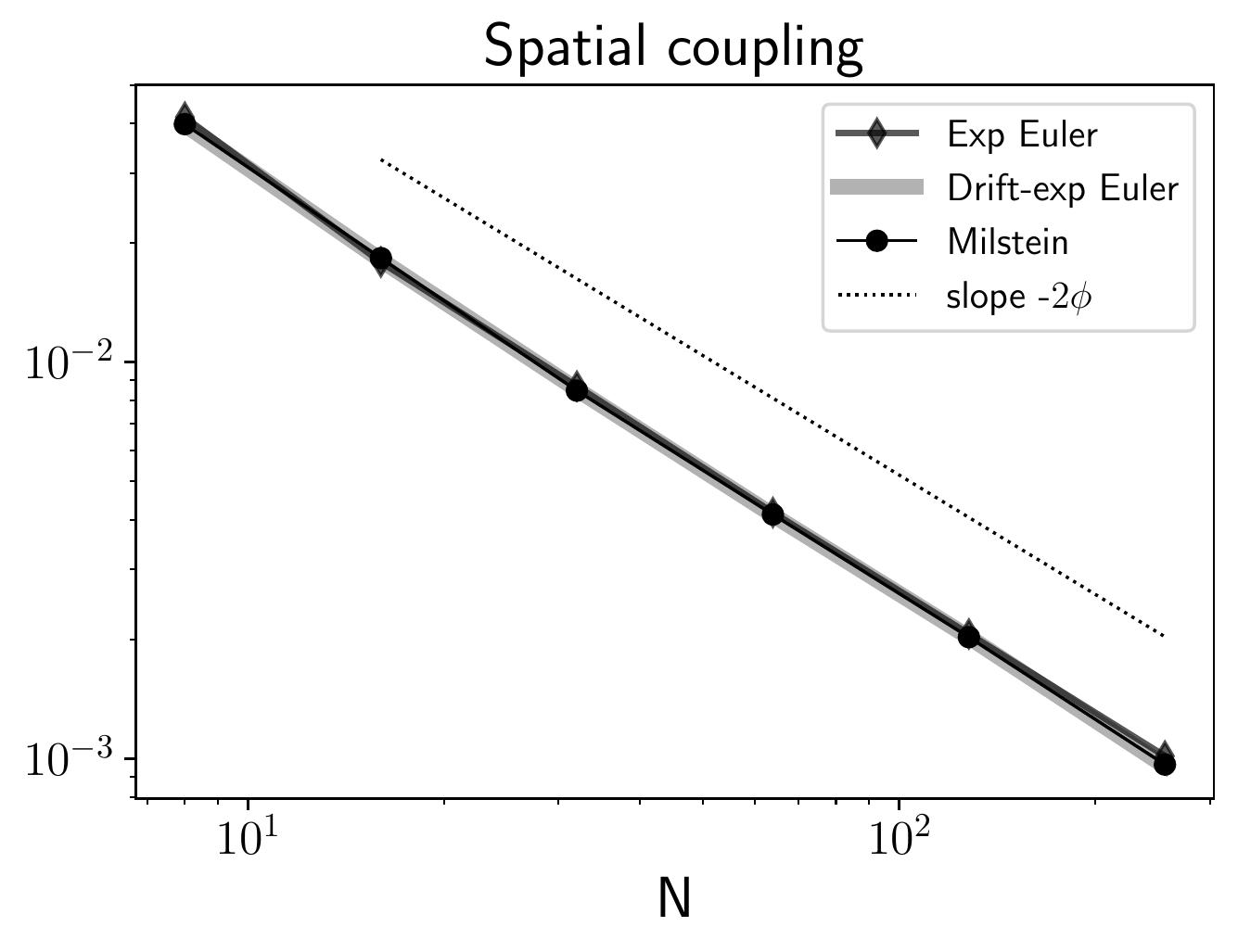}
  \includegraphics[width=0.49\textwidth]{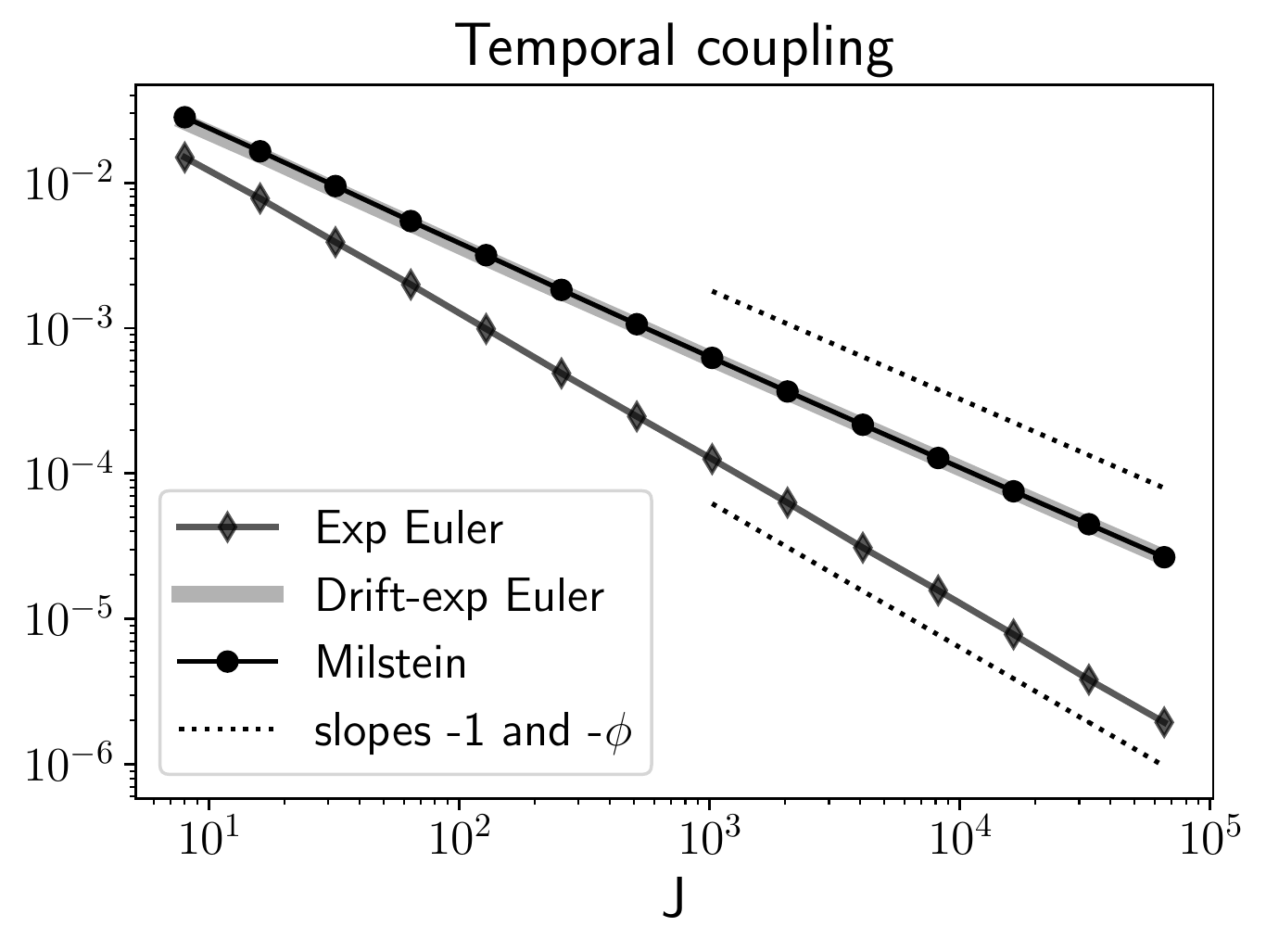}
  \includegraphics[width=0.49\textwidth]{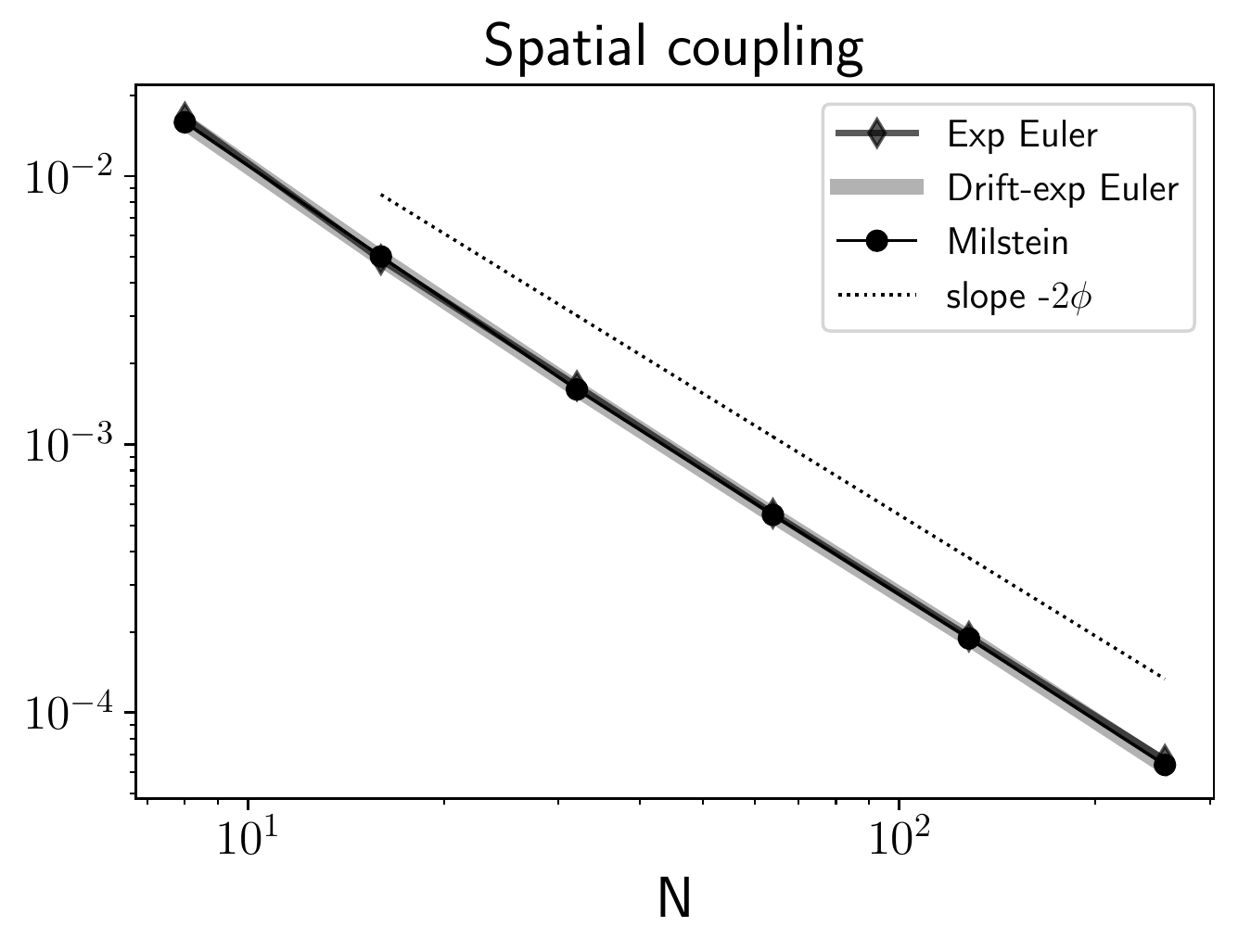}
  \caption{RMSE in time and space for the SPDE with the linear reaction term. The top row provides rates for low-regularity setting with $b=1/4$ and the bottom row provides rates for $b=1/2$. }
  \label{fig:convRatesLinear}
\end{figure}

\begin{figure}[h!] \center
  \includegraphics[width=0.49\textwidth]{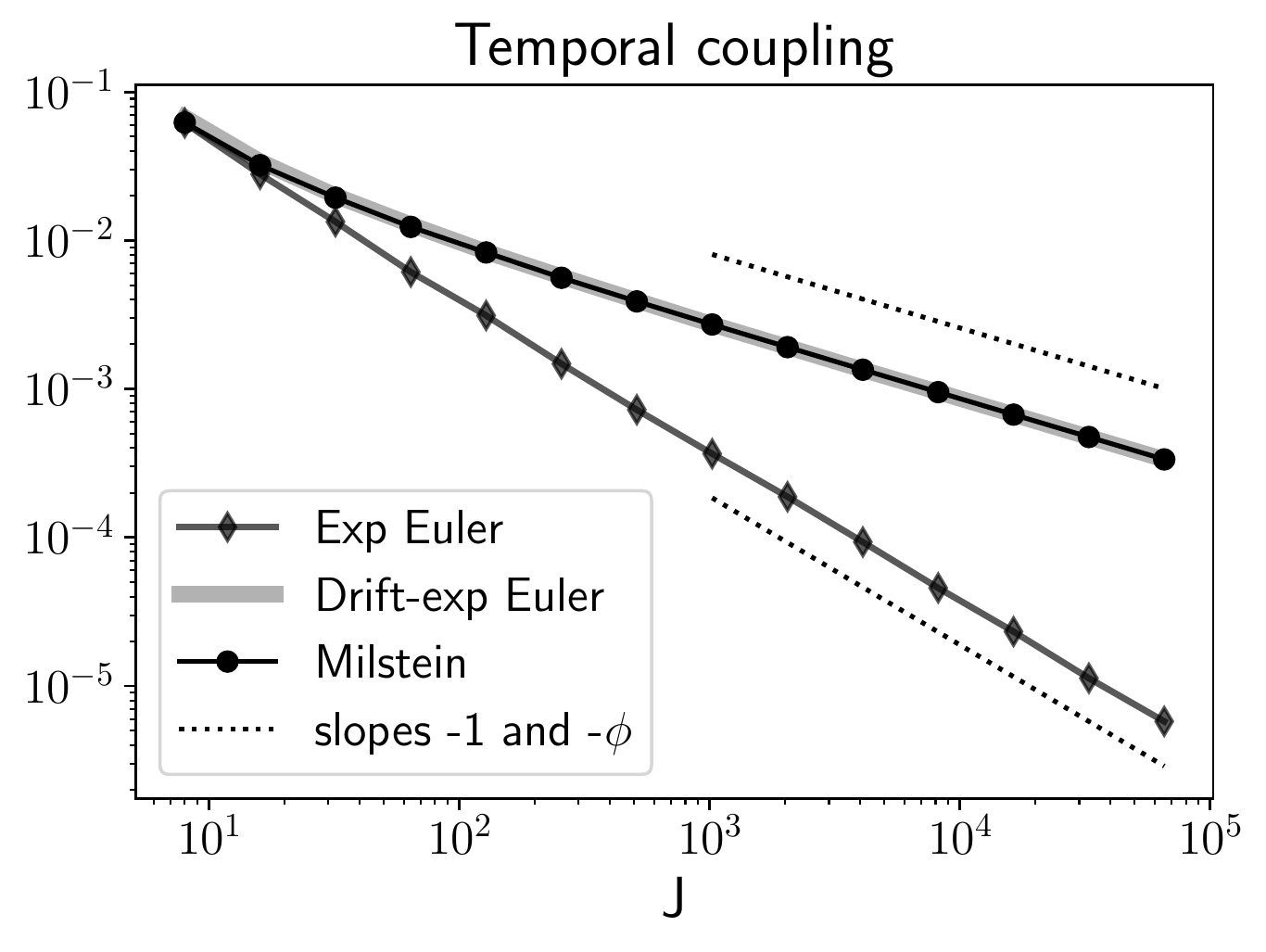}
  \includegraphics[width=0.49\textwidth]{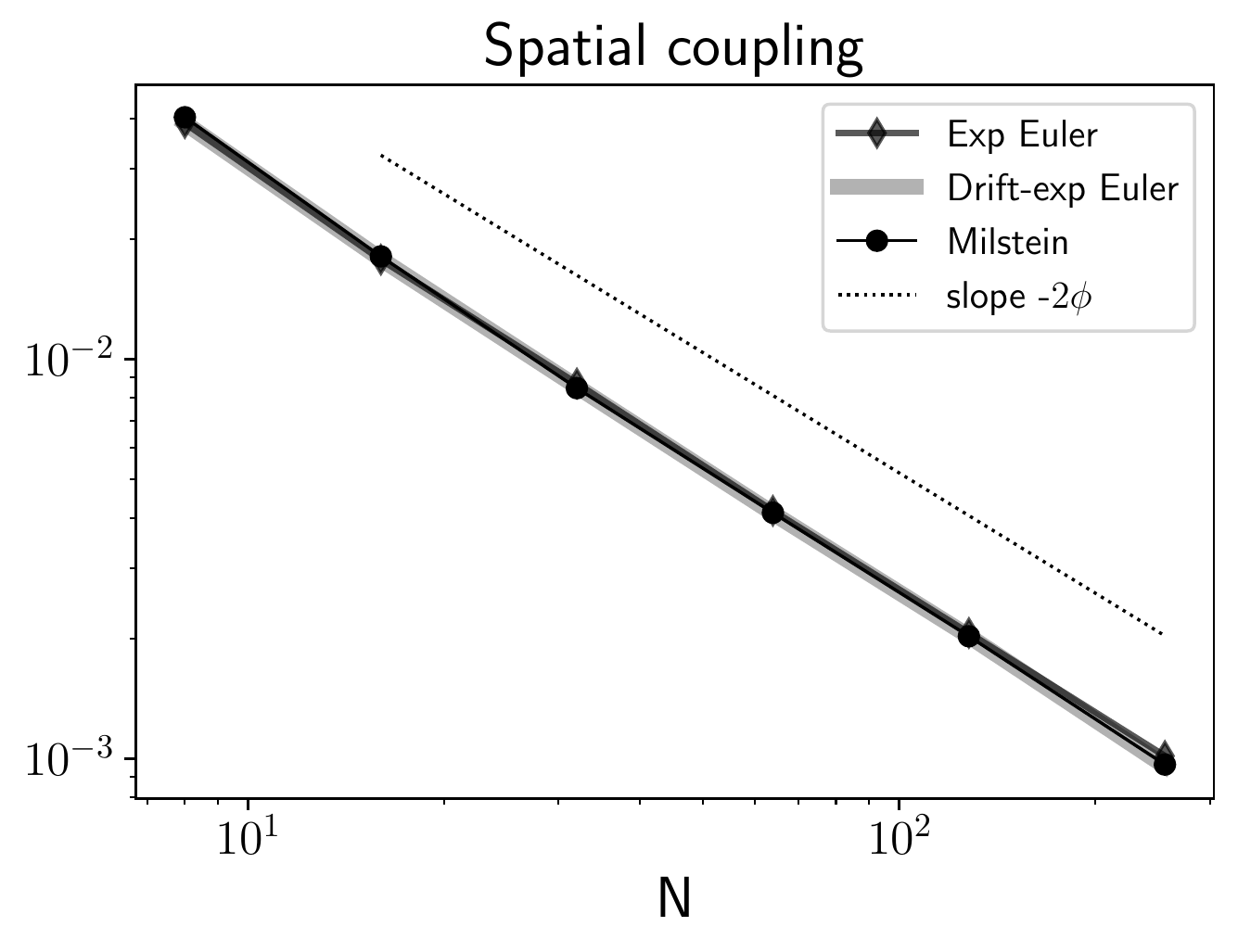}
  \includegraphics[width=0.49\textwidth]{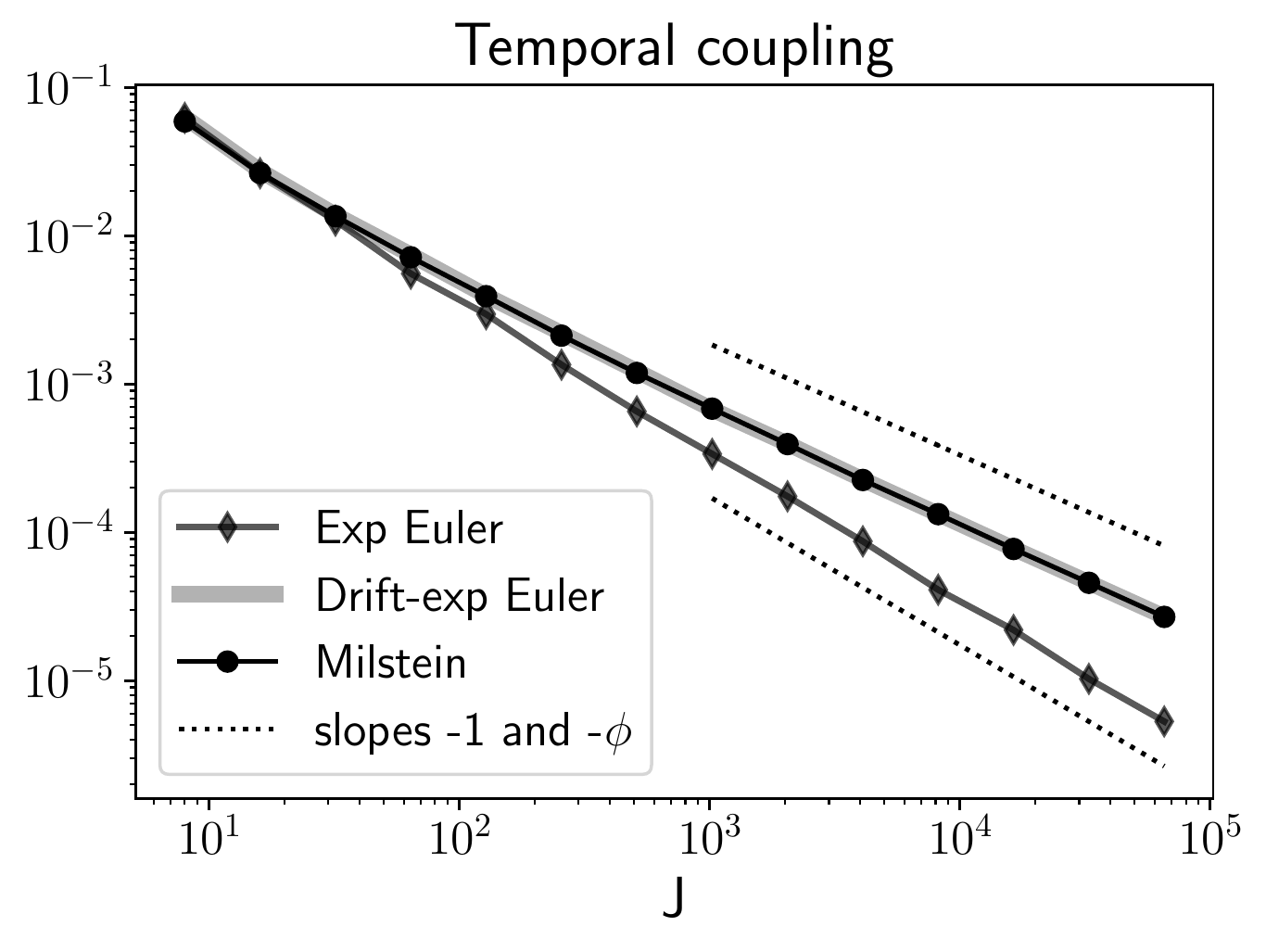}
  \includegraphics[width=0.49\textwidth]{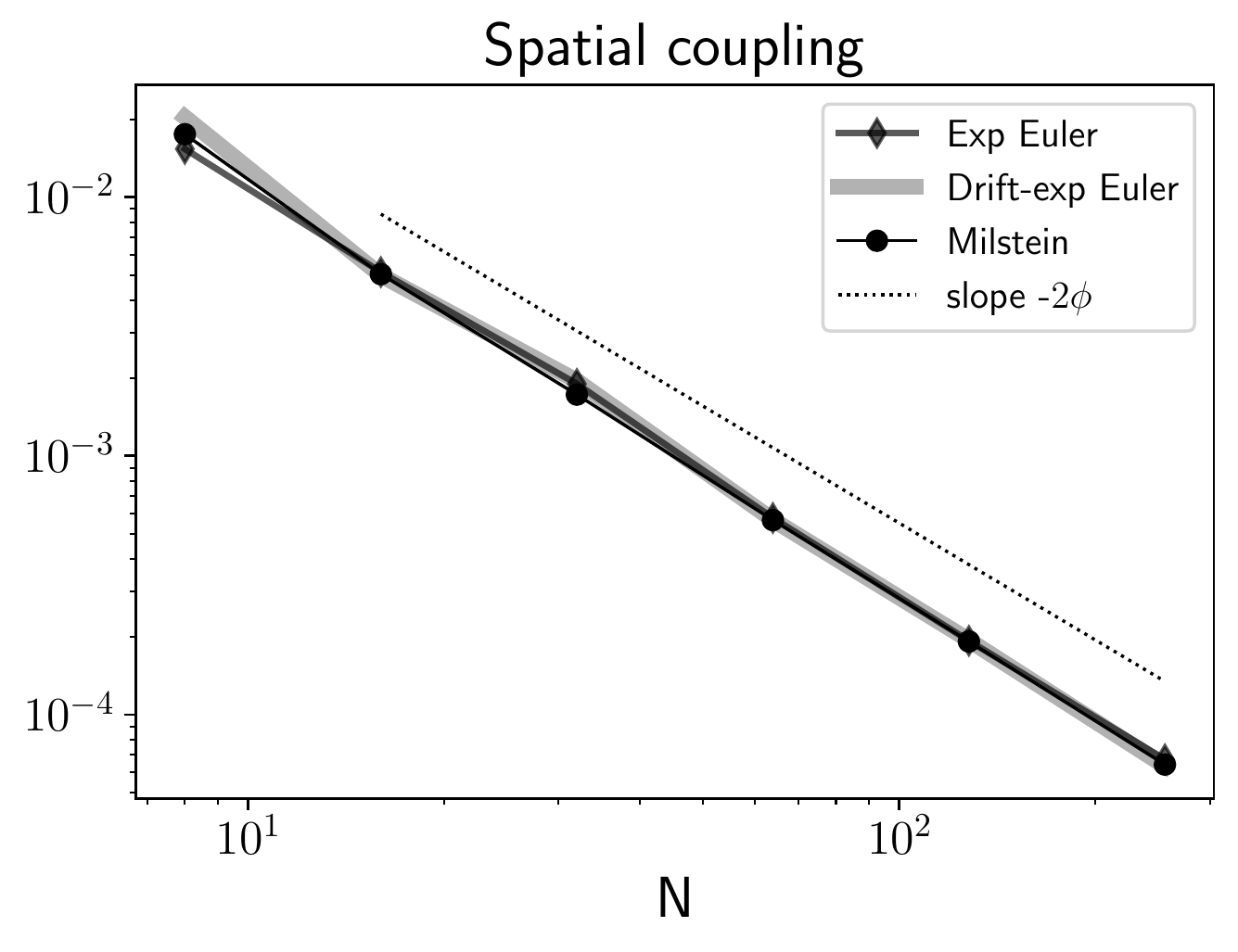}
  \caption{RMSE in time and space for the SPDE with the trigonometric reaction term. The top row provides the rates for low-regularity setting with $b=1/4$ and the bottom row provides the rates for $b=1/2$. }
  \label{fig:convRatesTrigonometric}
\end{figure}

Numerical estimates of the root mean squared error (RMSE) convergence rates
in time and space for all three methods are presented in Figures~\ref{fig:convRatesLinear} and~\ref{fig:convRatesTrigonometric}.
The RMSE in time is approximated by
\[
\sqrt{E_{M}[ \|V^{N_*, 2J}(T, \cdot) - V^{N_*, J}(T, \cdot) \|_H^2 ]}, 
\]
where $J$ is varied and $N_* = 1024$ is fixed, and using $M=10000$ independent samples of the random variable in the Monte Carlo estimator. For the exponential Euler method we observe the rate $1$ and for the other methods, we observe the rate $\phi(b) = 1/4+b$.

The RMSE in space is approximated by
\[
\sqrt{E_{M}[ \|V^{2N, J_*}(T, \cdot) - V^{N, J_*}(T, \cdot) \|_H^2 ] },
\]
where $N$ is varied and $J_* = 2^{18}$ is fixed, and using $M=250$ independent samples. This error describes the RMSE convergence rate in $N$, which we observe to be $2\phi = 1/2+2b$ for all methods.

Since the $\beta$ in Theorem~\ref{thm:VMLMC} represents the MSE, 
the numerical experiments indicate that $\beta = 2\min(1, 2\phi) = 2$ for exponential Euler MLMC and 
$\beta = 2\phi(b) = 1/2 + 2b$ for the other two methods. We will further set $\alpha = \beta/2$ as 
the weak rate when implementing all MLMC methods.

\subsection{Method parameters}

All three methods are implemented using Theorem~\ref{thm:VMLMC} 
with the numerical estimates of the rates $\alpha$ and $\beta$,
rather than by using the rather than using the slightly more 
conservative rate for $\beta$ in Theorem~\ref{thm:mainExpEuler} (ii). 

\subsubsection{Exponential Euler MLMC method}
We use the estimated rates $\alpha =1$ and $\beta=2$ and
balance error contributions in time and space by setting  
\[
  J_\ell = 2^{\ell+2} \quad \text{and} \quad N_\ell = 2 \times \lceil 2^{\ell/(2 \phi) +1} \rceil \quad  \text{and} \quad  .
\]
and $L = \lceil \log_2( 1/\epsilon )/\alpha \rceil = \lceil \log_2( 1/\epsilon )\rceil$.
We set $C_\ell := (\ell +1) 2^{\gamma \ell}$ with $\gamma = 1 + 1/(2\phi)$,
which one may verify is consistent with $C_\ell \eqsim J_\ell N_\ell \log_2(N_\ell)$, and we set $V_\ell := 2^{-2\ell}$ 
to determine the sequence $\{M_\ell\}_\ell$, in compliance with formula~\eqref{eq:ml2}, by 
\begin{equation} \label{eq:ml_values}
M_\ell(\epsilon) = \begin{cases}
  20 \left\lceil 2\epsilon^{-2} \sqrt{\frac{V_\ell}{C_\ell}}
    \sum_{j=0}^L\sqrt{V_j C_j} \right\rceil & \text{if} \quad \ell = 0\\[5pt]
    \; \;5 \left\lceil 2\epsilon^{-2} \sqrt{\frac{V_\ell}{C_\ell}}
    \sum_{j=0}^L\sqrt{V_j C_j} \right\rceil & \text{otherwise.}
\end{cases} 
\end{equation}

\subsubsection{Drift-exponential Euler MLMC and Milstein MLMC }
The numerically observed convergence rates for
both of these methods are $\beta = 2\phi$ and $\alpha = \phi$.  
For both methods, we set 
\[
  J_\ell = 2^{\ell+2}, \quad N_\ell = 2\times \lceil 2^{\ell/2 +1}\rceil,
  \quad \text{and} \quad L = \left\lceil \frac{\log_2(1/\epsilon)}{\phi} \right\rceil-2,
\]
$C_\ell := (\ell +1) 2^{3\ell/2}$ and $V_\ell := 2^{- 2\phi \ell}$,
we and determine the sequence $\{M_\ell\}_{\ell}$ by formula~\eqref{eq:ml_values}.

\subsection{Linear reaction term}

\begin{figure}[h!] \centering
  \begin{subfigure}[b]{0.49\textwidth}
  \includegraphics[width=\textwidth]{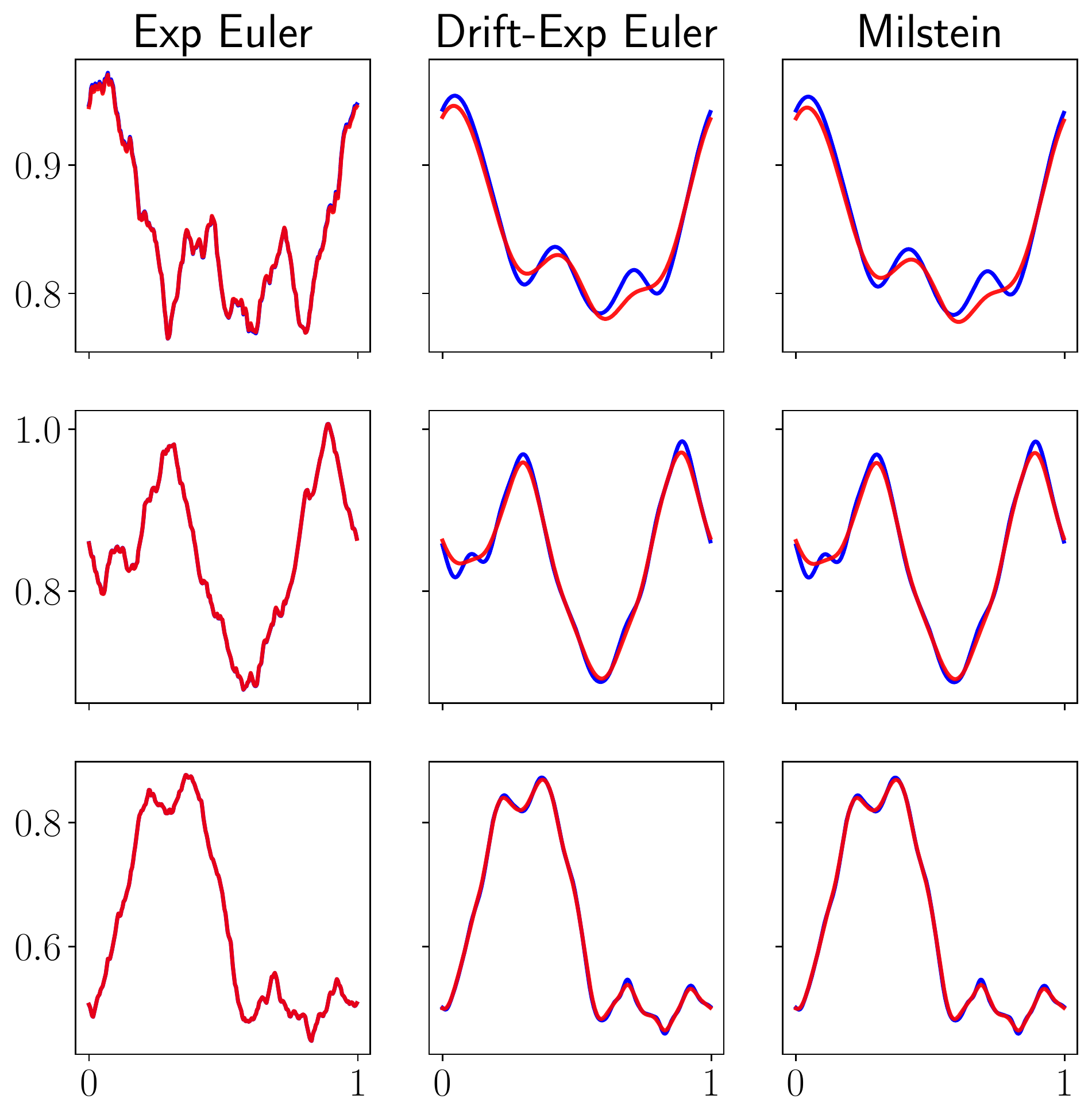}
    \caption{$b=1/4$}
\end{subfigure}
\begin{subfigure}[b]{0.49\textwidth}
  \includegraphics[width=\textwidth]{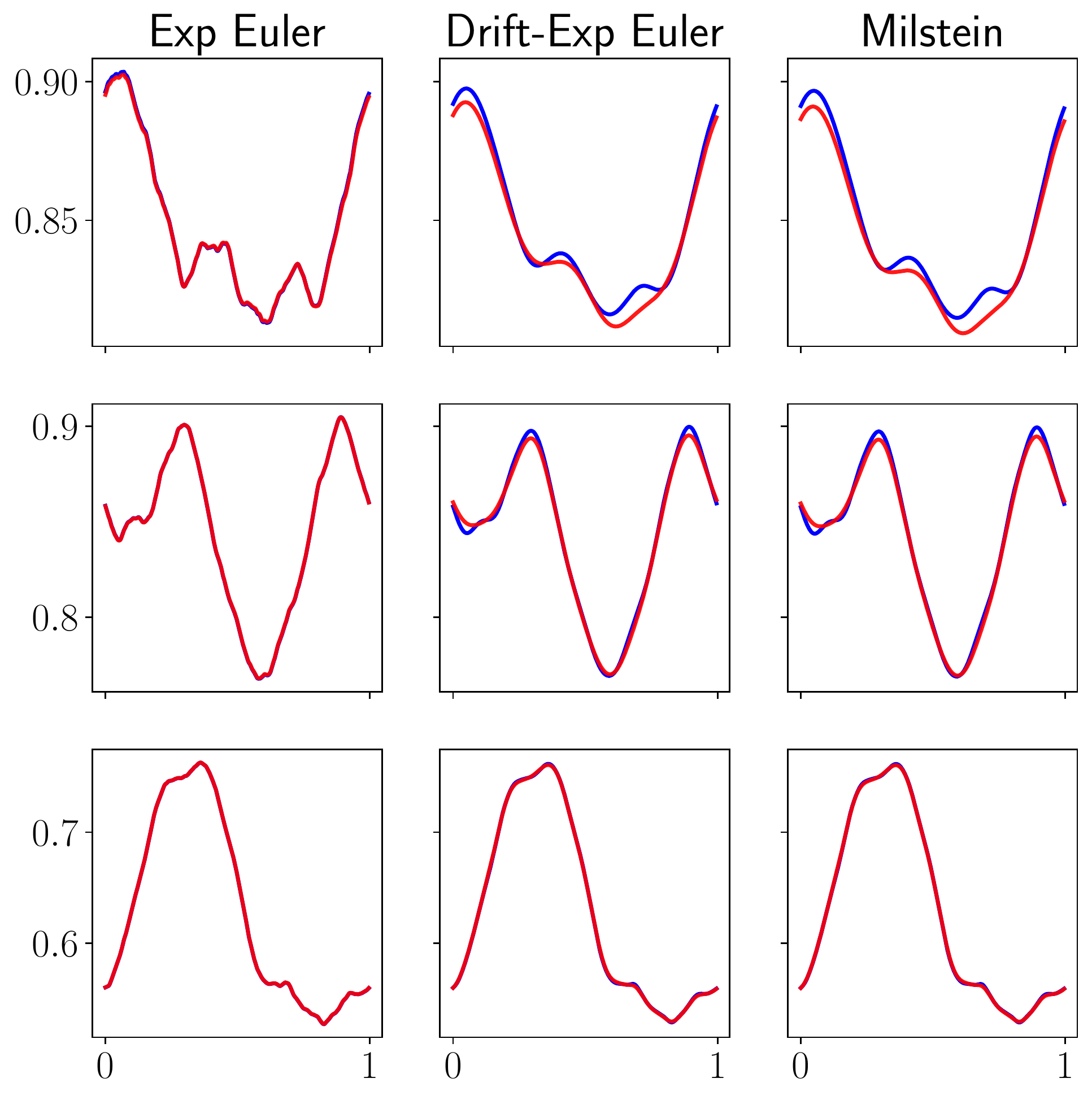}
  \caption{$b=1/2$}
\end{subfigure}
  \caption{Comparison of the couplings at the finest level (blue) and
    the coarsest level (red) for the SPDE with linear reaction term.
    The spatial resolution is fixed to $N=2^8$ in all simulations and the
    resolution in time $(J, J/2)$ is given by $J=2^6, 2^8,$ and $2^{10}$ from
    top to bottom row.}
  \label{fig:linb025Coupling}
\end{figure}

We first consider the SPDE with $f(U) = U$. 
Figure~\ref{fig:linb025Coupling}
presents pairwisely coupled realizations for 
progressively finer resolution in time for 
the settings $b=1/4$ and $b=1/2$. 
In the low regularity setting $b=1/4$, we clearly observe that the exponential 
Euler method has far less smoothing of the solutions and achieves a stronger coupling
than the other methods. The difference between the methods becomes less 
visible in the smoother setting $b=1/2$.

\begin{figure}[h!] \centering
  \includegraphics[width=0.75\textwidth]{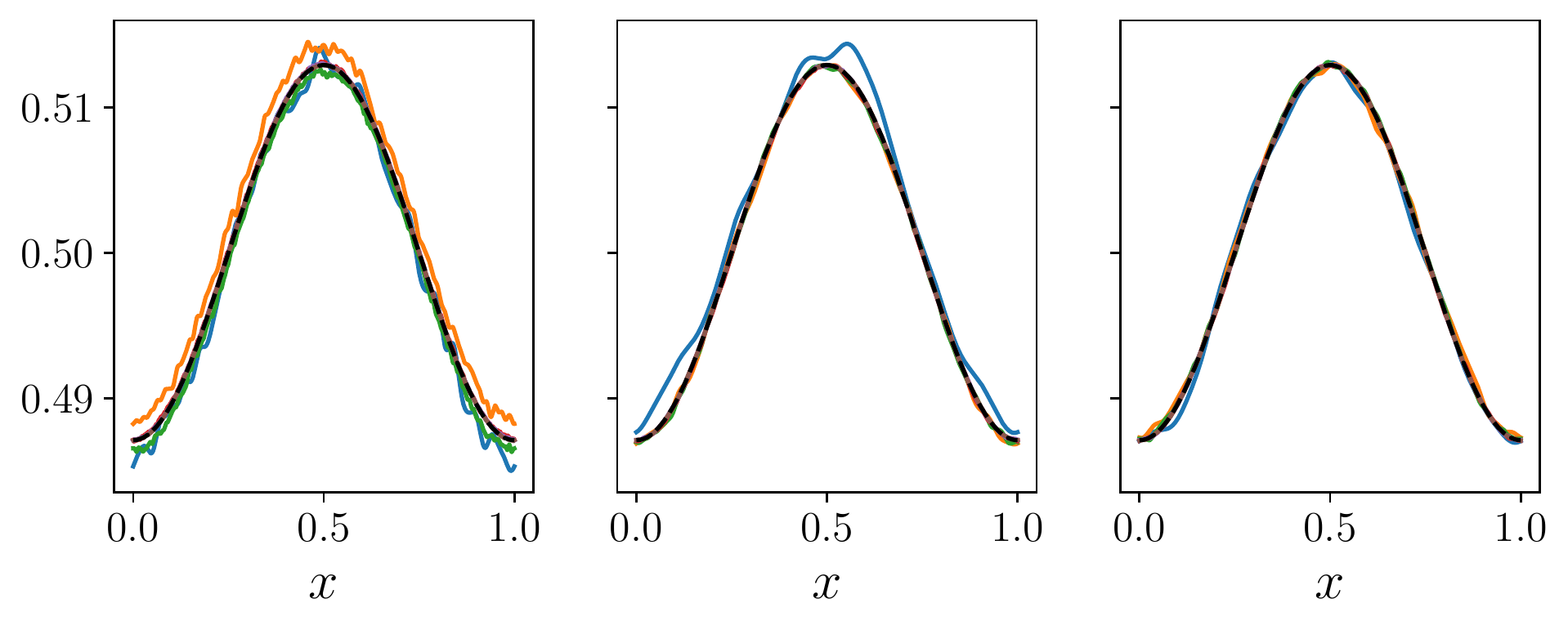}
  \includegraphics[width=0.75\textwidth]{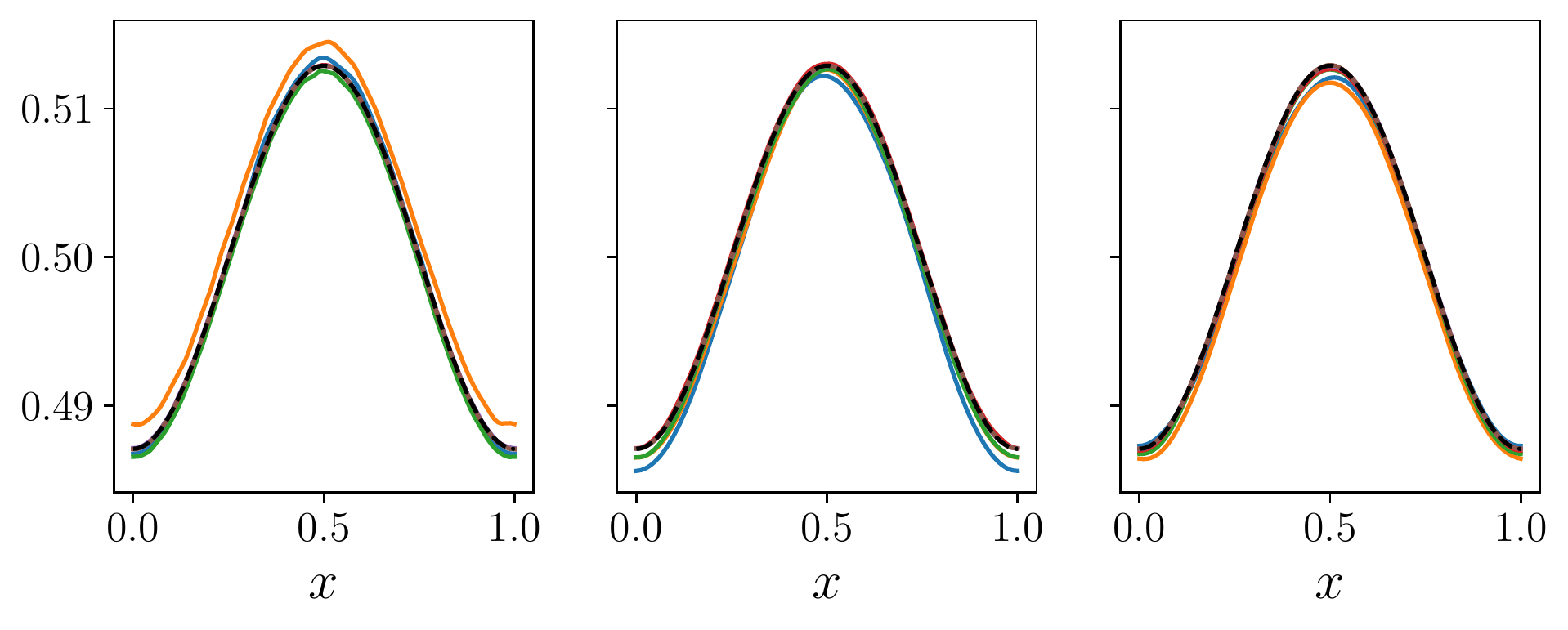}
   \caption{Approximation of $E_{\mathrm{ML}}[U(T,x)]$ 
    for the SPDE with linear reaction term for different
    values of $\epsilon = 2^{-4}, 2^{-5}, \ldots, 2^{-9}$ in full
    colored lines (blue, orange, green, red, purple, brown)
    and the pseudo-reference solution $\mathbb{E}[U(T,x)]$ (dashed line).  
    From left to right, exponential Euler, drift-exponential Euler, and Milstein.
    Top row is for the low-regularity setting $b=1/4$ and bottom row is for 
    $b=1/2$.    
    }
  \label{fig:linb025ApproxMLMC}
\end{figure}

Figure~\ref{fig:linb025ApproxMLMC} provides the approximation of
$\bE[U(T,\cdot)]$ by one simulation of each of the MLMC methods 
for different input $\epsilon = 2^{-\ell}$ for $\ell = 4,5,\ldots,9$. We observe that 
all methods converge to the mean with approximately the same rate
for both values of $b$.  

\begin{figure}[h!] \center
  \includegraphics[width=0.49\textwidth]{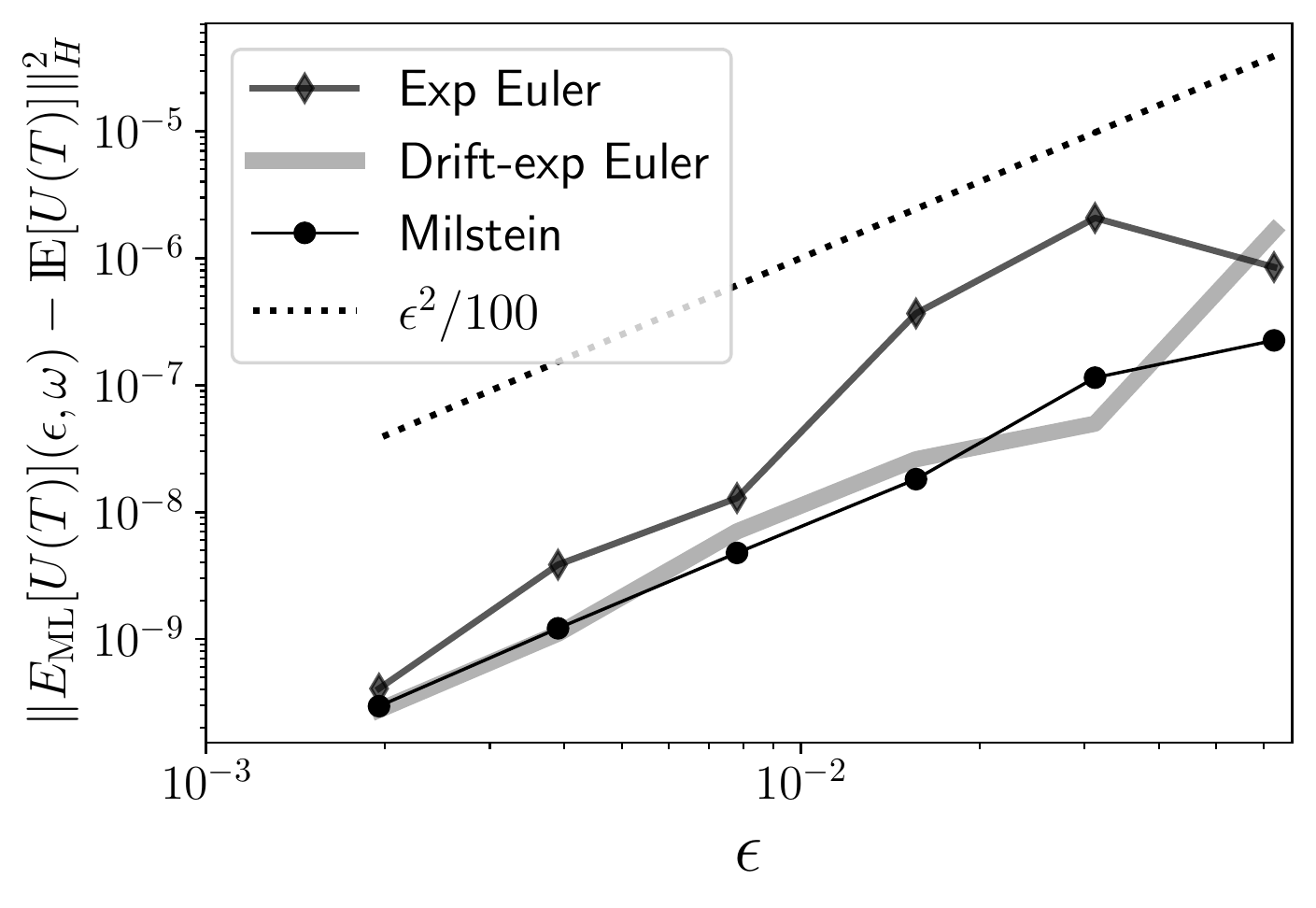}
  \includegraphics[width=0.49\textwidth]{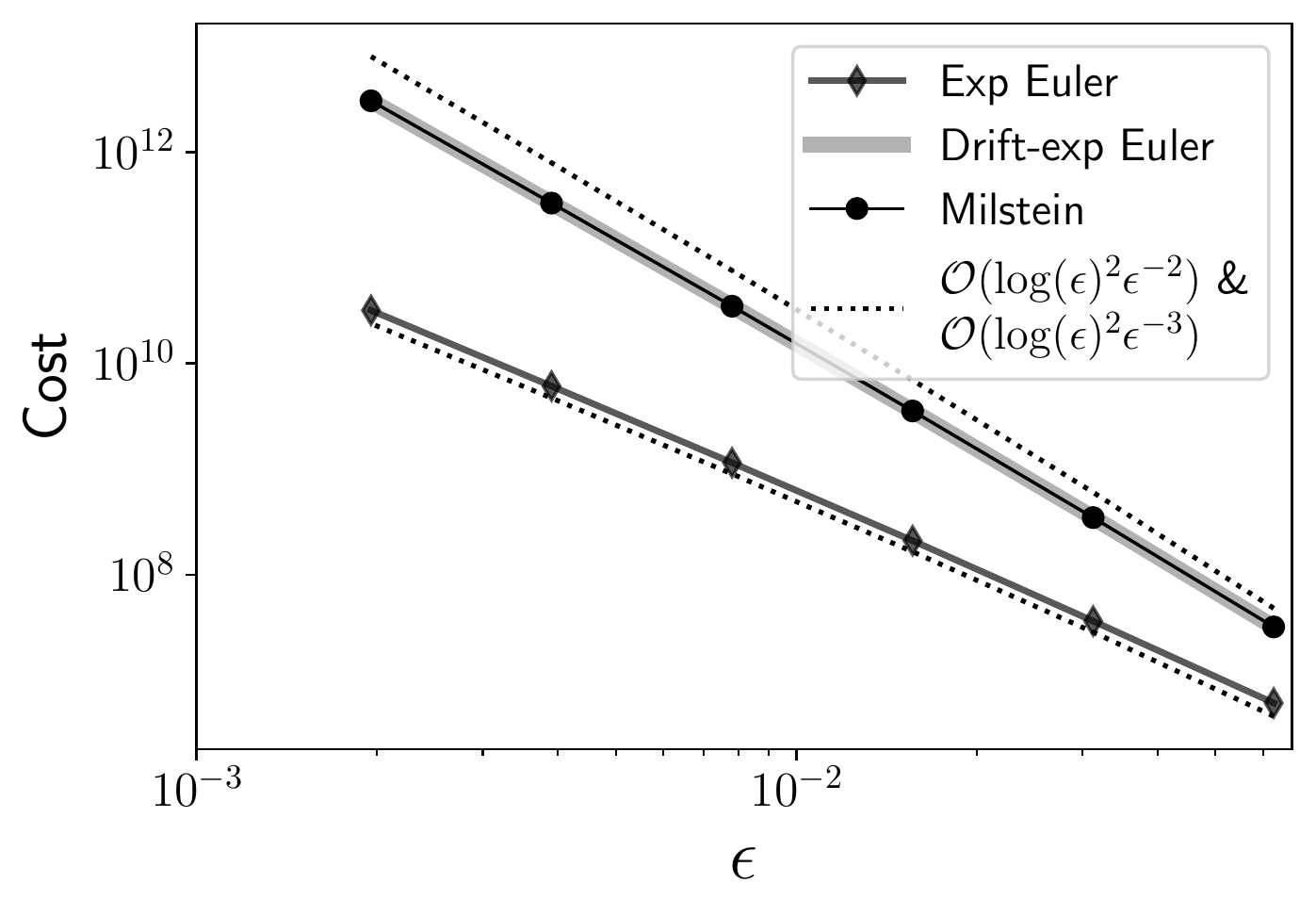}
  \includegraphics[width=0.49\textwidth]{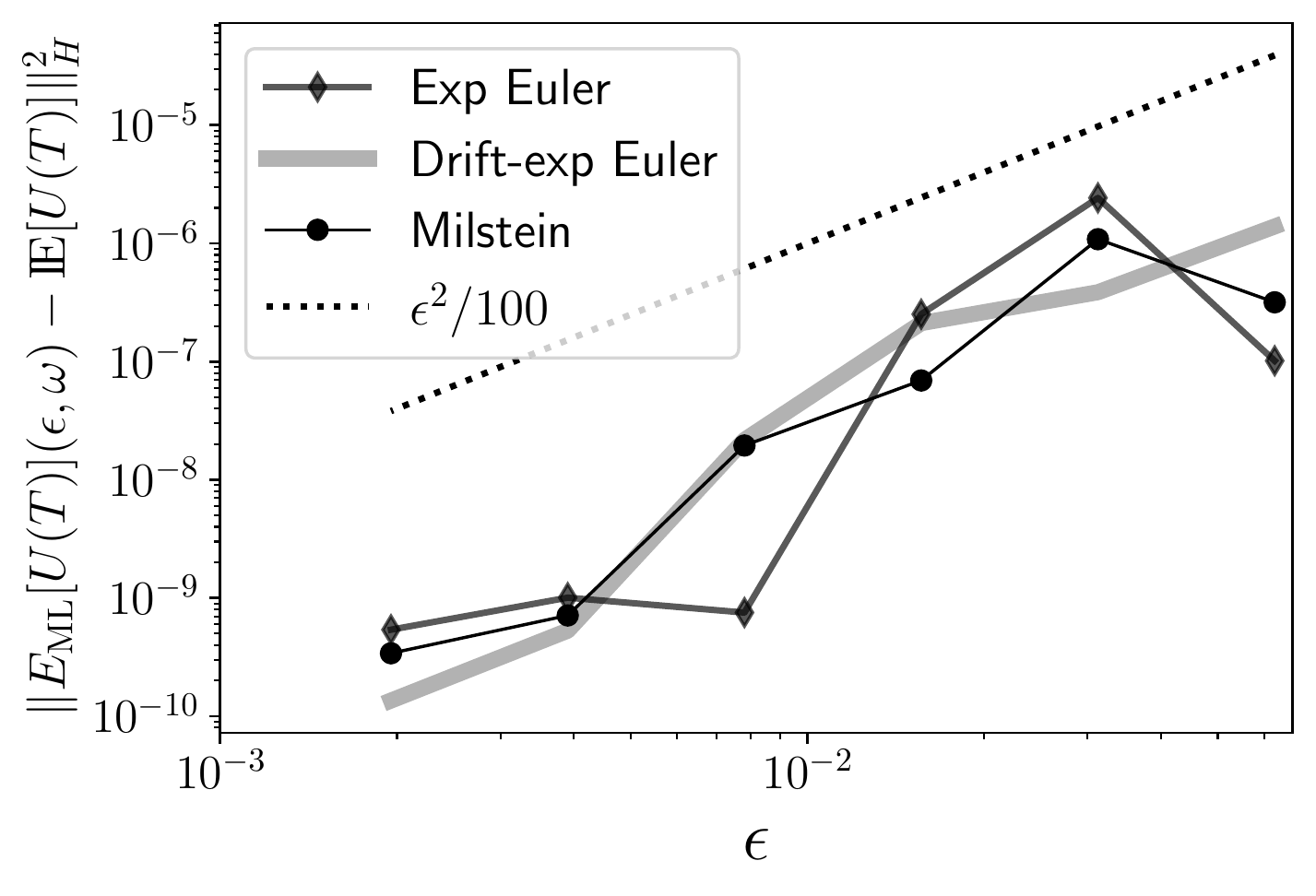}
  \includegraphics[width=0.49\textwidth]{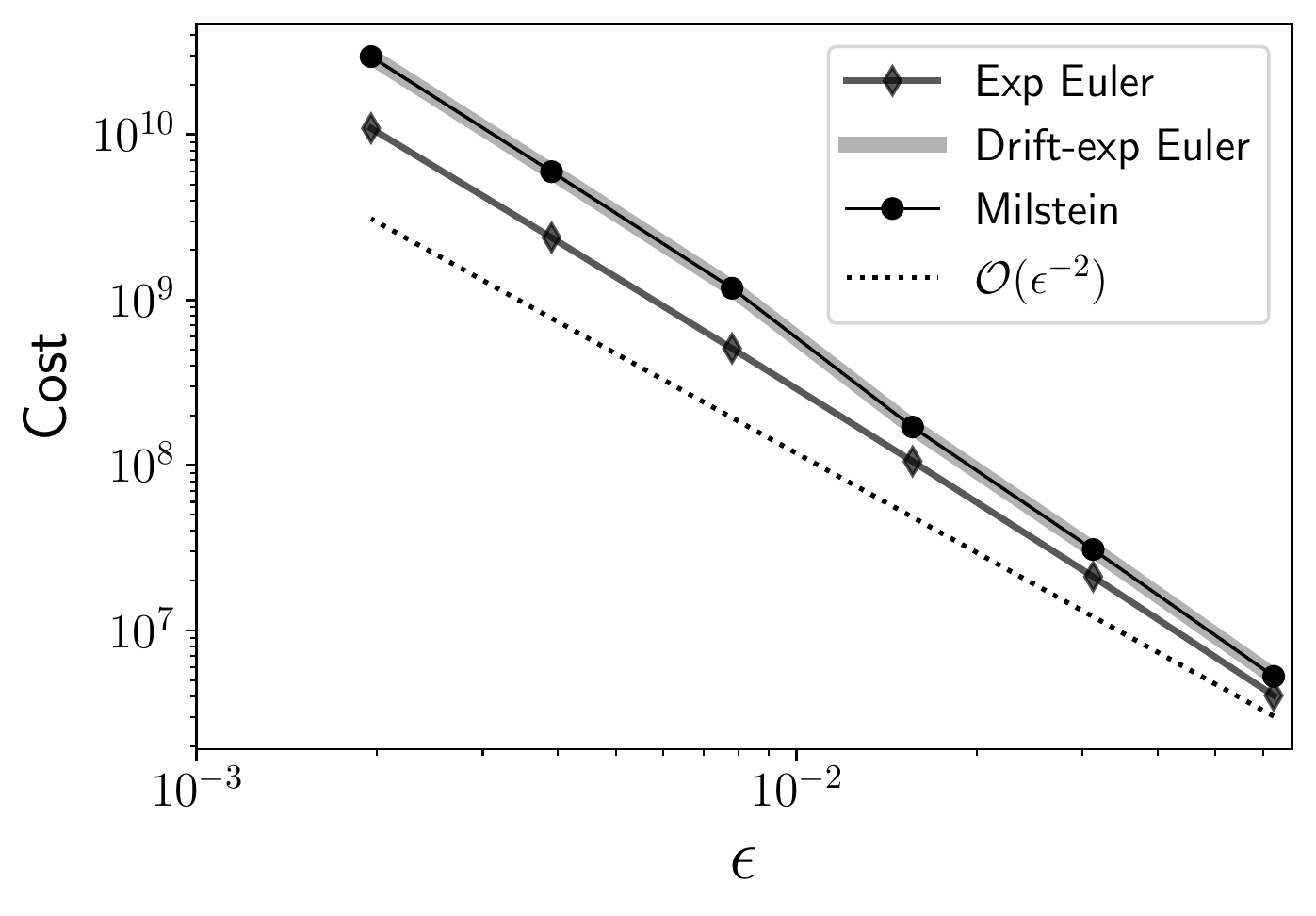}
  \caption{Top row: convergence and computational cost plots for the SPDE with 
  $f(U)=U$ and $b=1/4$. Bottom row: similar plots for the setting with $b=1/2$.}
  \label{fig:linb025ErrorMLMC}
\end{figure}

Figure~\ref{fig:linb025ErrorMLMC} presents the MLMC approximation
error versus tolerance and the computational cost versus tolerance for
different input tolerances $\epsilon$ for one simulation. 
The approximation error 
\[
  \|E_{\rm{ML}}[U(T,\cdot)](\omega; \epsilon) - \E{U(T, \cdot)} \|_H^2
\]
is computed for one simulation of the MLMC estimator for each input of $\epsilon$, where
the the pseudo-reference solution $\E{U(T, \cdot)}$ is obtained by solving the PDE 
\[
  dU_t^N = (A_N U_t^N + f_N(U^N)) dt, \qquad U_0^N = P_N u_0
\]
with the exponential Euler method using the resolutions $N=2^{13}$ and $J=2^{18}$. Let us also recall that the computational cost of the MLMC methods is defined by $\sum_{\ell=0}^L C_\ell M_\ell$.  

For $\phi = 1/2-$, we observe that exponential Euler MLMC method has achieves the error $\mathcal{O}(\epsilon^2)$ 
at the cost $\mathcal{O}((\log_2(\epsilon))^2 \epsilon^{-2})$ while the other methods achieves similar accuracy at  considerably higher cost. For $\phi = 3/4-$ all three methods achieves an error $\cO(\epsilon^2)$ at
a comparable computational cost. The observations are consistent with theory. 

\subsection{Trigonometric reaction term}
We next consider the SPDE with 
\[
f(U)(x) = 2\big(\sin(2 \pi U(x)) + \cos(2 \pi U(x)) \big).
\]

\begin{figure}[h!] \centering
  \includegraphics[width=0.75\textwidth]{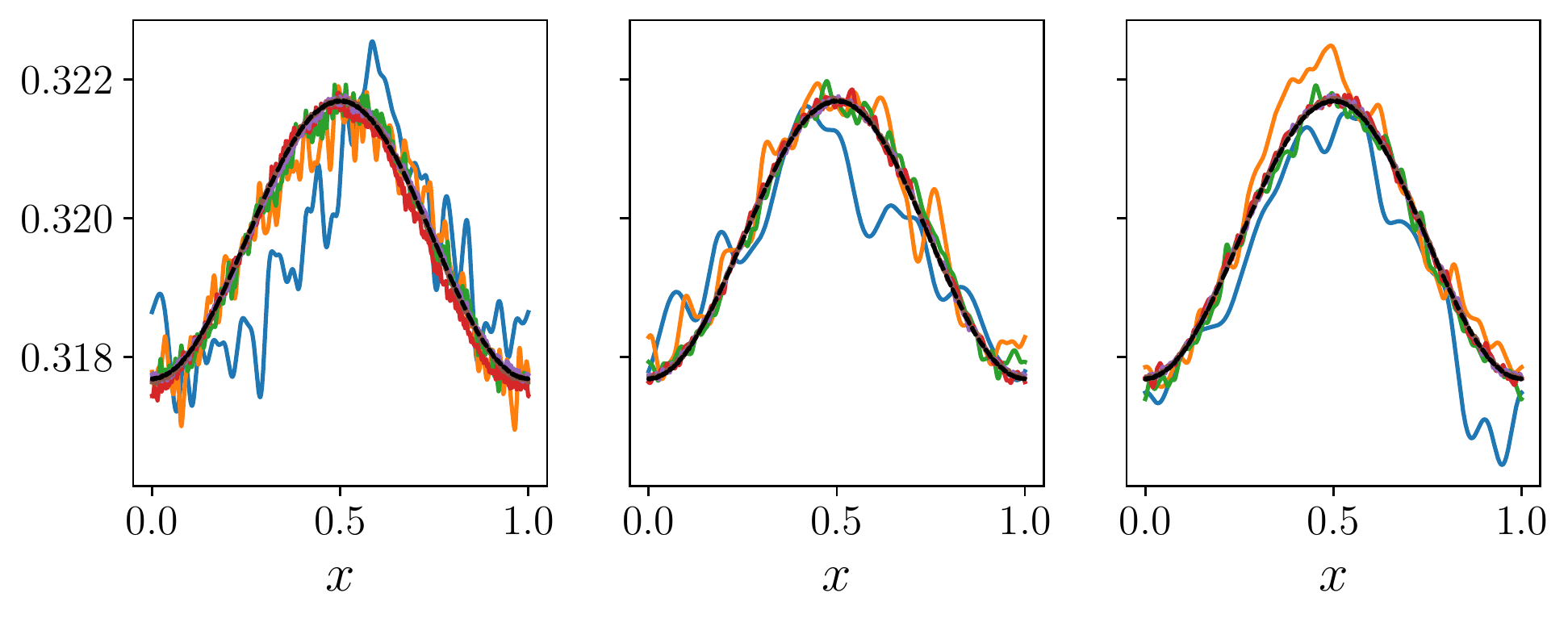}
  \includegraphics[width=0.75\textwidth]{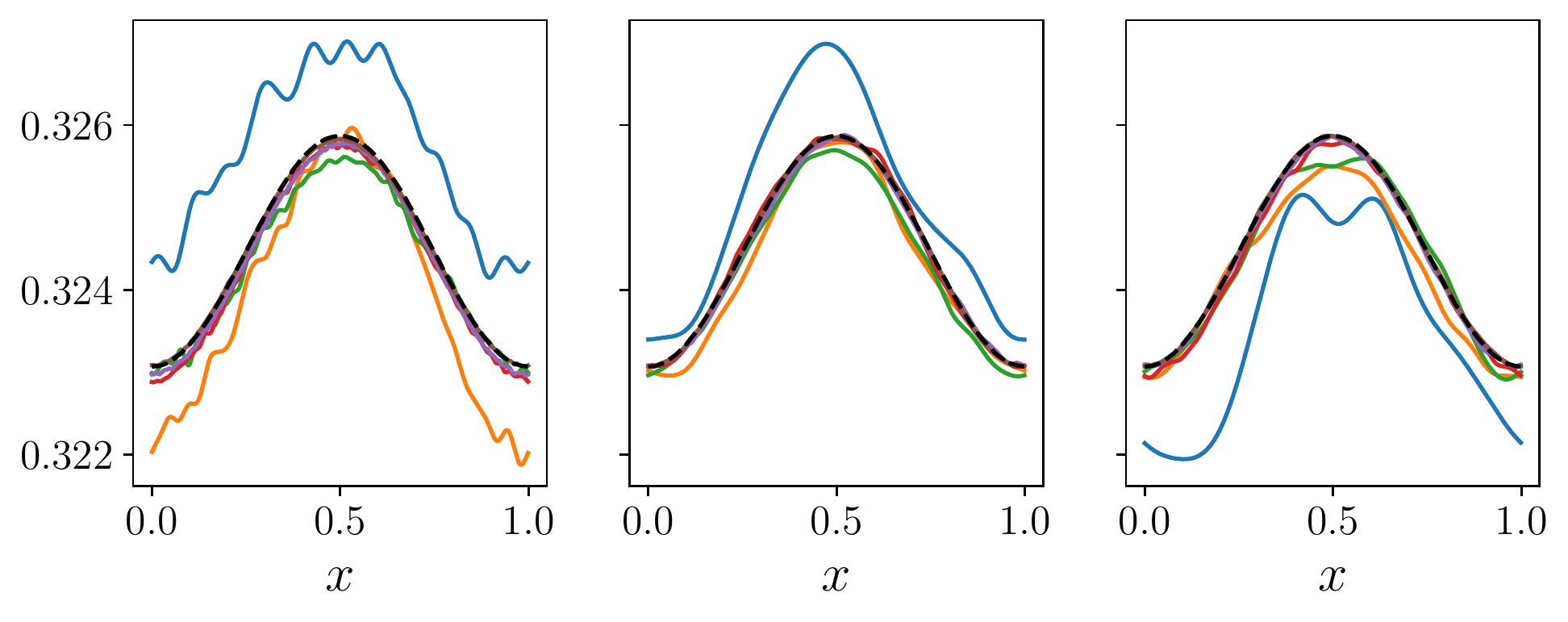}
  \caption{Approximation of $E_{\mathrm{ML}}[U(T,x)]$ 
    for the SPDE with trigonometric reaction term for different
    values of $\epsilon = 2^{-4}, 2^{-5}, \ldots, 2^{-9}$ in full
    colored lines (blue, orange, green, red, purple, brown)
    and the pseudo-reference solution $\mathbb{E}[U(T,x)]$ (dashed line).  
    From left to right, exponential Euler, drift-exponential Euler, and Milstein.
    Top row is for the low-regularity setting $b=1/4$ and bottom row is for 
    $b=1/2$.    
    }
  \label{fig:sinb025ApproxMLMC}
\end{figure}

\begin{figure}[h!] \centering
  \includegraphics[width=0.49\textwidth]{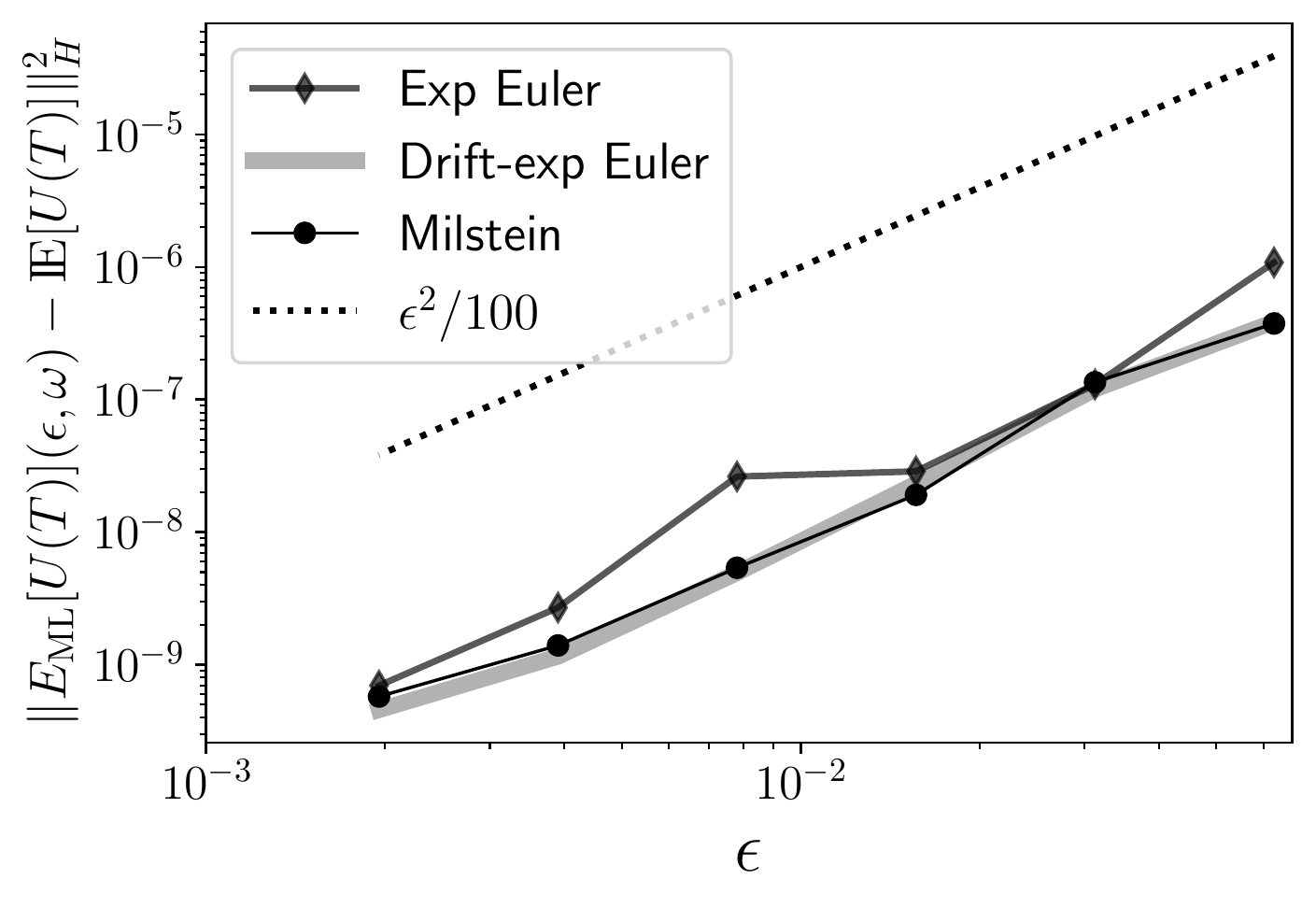}
  \includegraphics[width=0.49\textwidth]{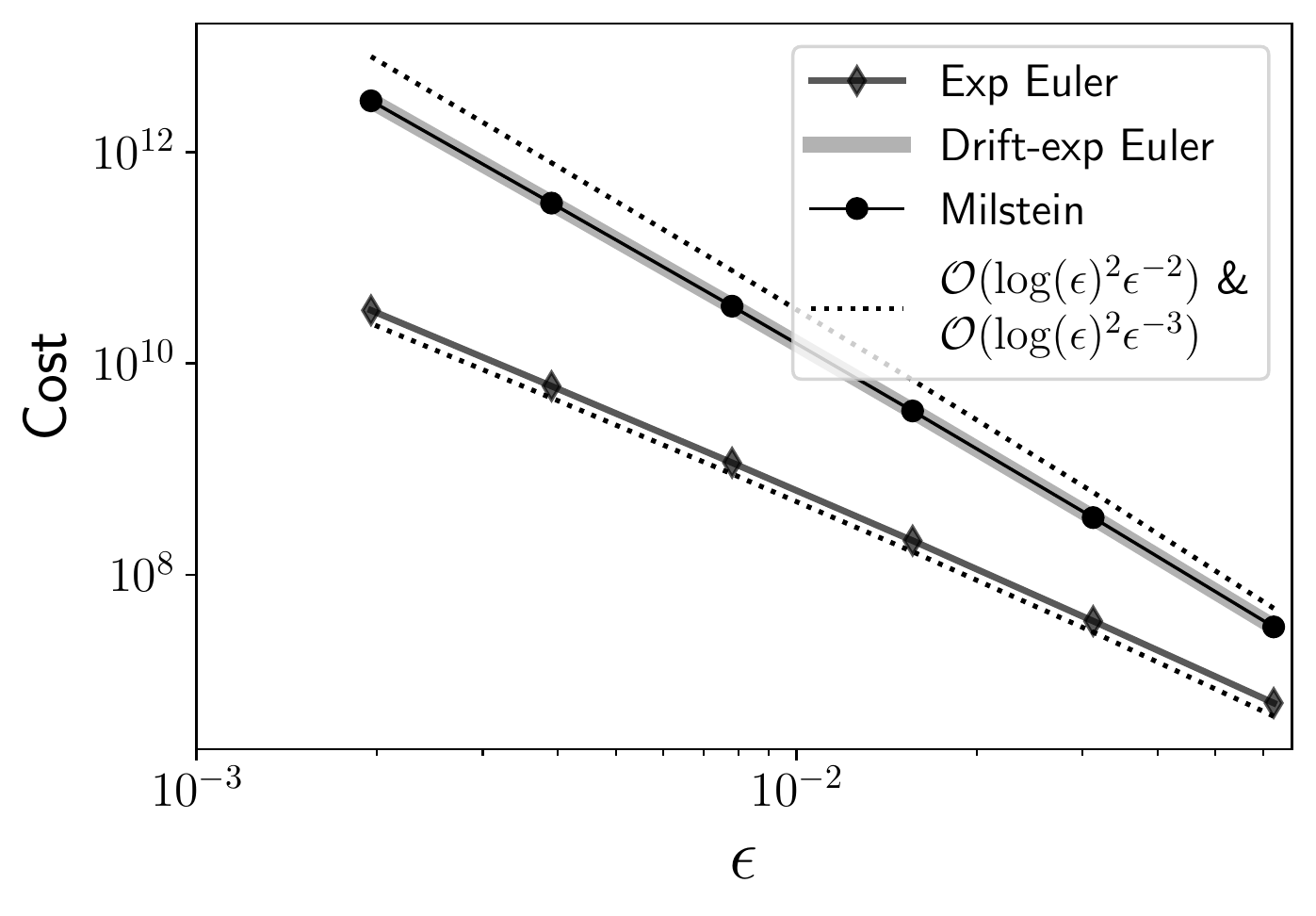}
  \includegraphics[width=0.49\textwidth]{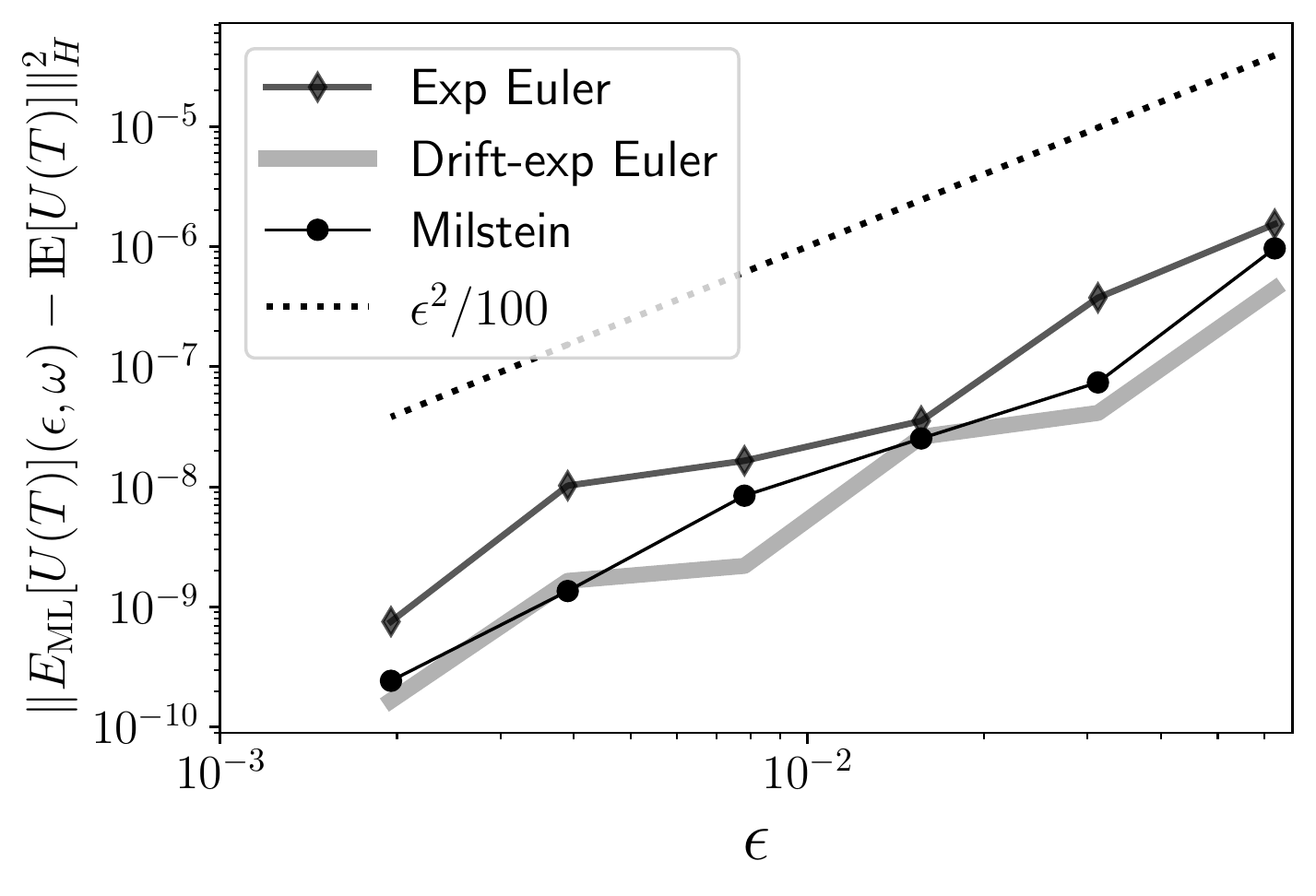}
  \includegraphics[width=0.49\textwidth]{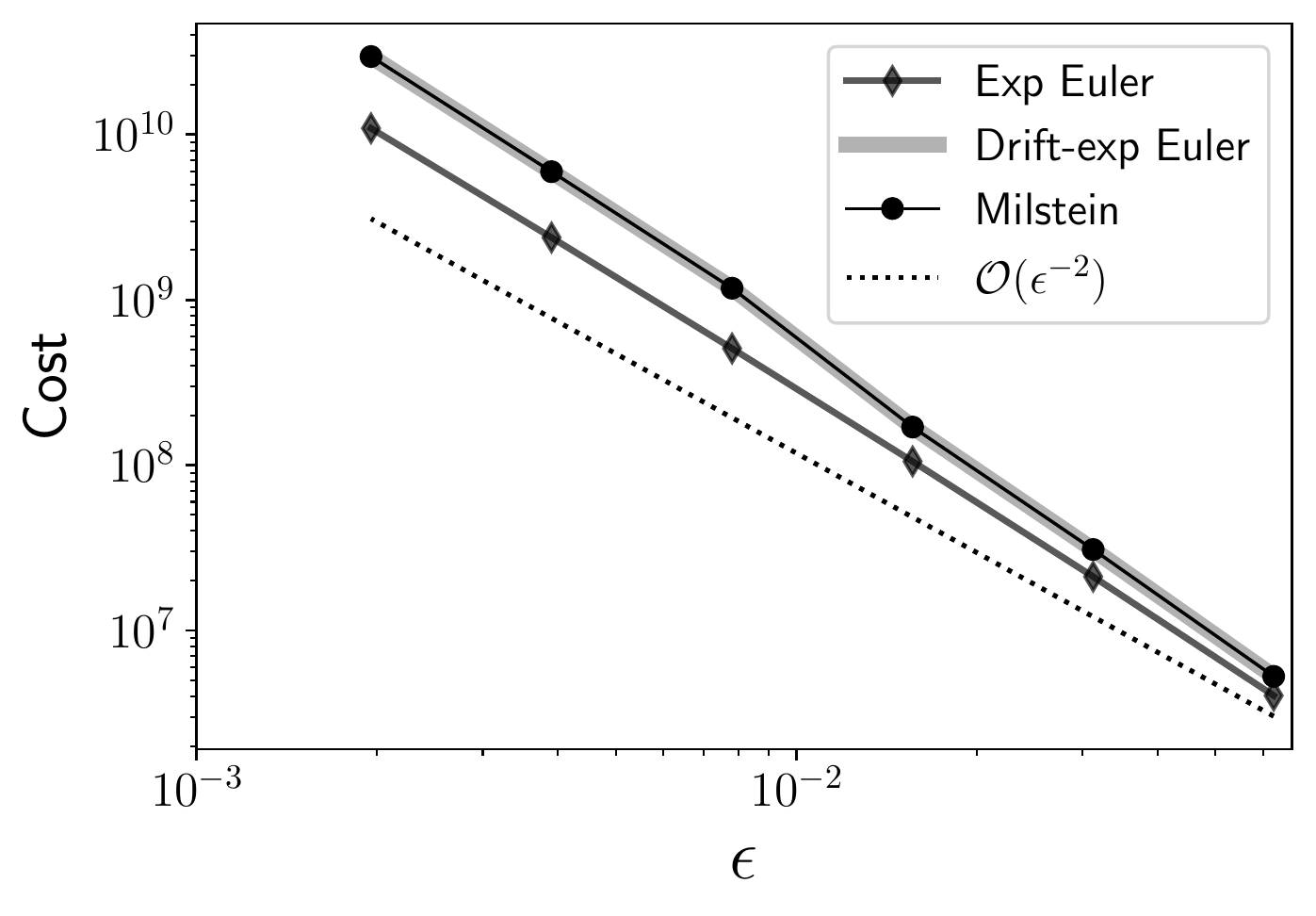}
  \caption{Top row: convergence and computational cost plots for the SPDE with 
  $f(U)= 2(\sin(2\pi U) + \cos(2\pi U)) $ and $b=1/4$. 
  Bottom row: similar plots for the setting with $b=1/2$.}
  \label{fig:sinb025ErrorMLMC}
\end{figure}

The approximation of $\bE[U(T,\cdot)]$ by the MLMC methods for different inputs
$\epsilon = 2^{-\ell}$ for $\ell = 4,5,\dots,9$ is presented in Figure~\ref{fig:sinb025ApproxMLMC}
and Figure~\ref{fig:sinb025ErrorMLMC} shows the MLMC approximation error 
versus computational cost and computational cost versus tolerance for
different input tolerances $\epsilon$. 
For each value of $b$, the pseudo-reference solution used for evaluating the 
approximation error is computed by the exponential Euler MLMC method
$E_{\mathrm{ML}}[ U(T,\cdot)](\omega, \epsilon) \approx \bE[U(T,\cdot)]$ 
with the overkilled parameter value $\epsilon = 2^{-11}$. 
This an expensive computation using the following number of samples per level when $b=1/4$:
\[
(M_0, M_1, M_2, \ldots, M_{10}, M_{11}) =   (4907168680, 216868270, 44268050, \ldots,  355, 85),
\]
with $N_{\ell} = J_{\ell} = 2^{\ell+2}$. We observe once again that
exponential Euler MLMC outperforms the other methods in the 
low-regularity setting $b=1/4$ and that all methods perform 
similarly when $b=1/2$.

\section{Conclusion}
\label{sec:conc}

Our objective in this work was to show both theoretically and
experimentally that coupling approaches that exploit more information 
than only the driving noise $W_t$, such as the exponential Euler MLMC method, 
can result in strong coupling and improve the 
efficiency of weak approximations for SPDE. Our motivation in doing so, was based 
on the lack of literature on strong coupling for MLMC methods solving SPDE. 
In particular, we have derived explicit convergence rates, related to the decay of the mean squared
error-to-cost rate, for the exponential Euler MLMC 
method and the Milstein MLMC method, cf.~Theorems~\ref{thm:mainExpEuler} and~\ref{thm:milsteinMLMC}. 
The convergence rates for exponential Euler MLMC method 
is an improvement over existing MLMC methods for   
reaction-diffusion SPDE with additive noise. We also presented 
numerical experiments highlighting our derived rates and demonstrating the
efficiency gains of the exponential Euler 
MLMC method over alternative ones. 
This was tested numerically on SPDE with linear and nonlinear reaction 
terms.

There are many possible extensions of this work.  It would be
interesting to understand whether strong couplings also can improve
the efficiency of MLMC for other numerical solvers for SPDE, such as
finite difference methods and FEM \cite{ABS13,BLS13}. This indeed is
a challenging problem, due to the seemingly limitless
possibilities of couplings for infinite-dimensional problems.  Another
direction is to develop a multi-index
Monte Carlo method~\cite{HNT16,JLX21} based on the pathwise correctly
coupled exponential Euler method. This has the potential of further
improving tractability in higher-dimensional physical space and
low-regularity settings.

 
\appendix
\section{Model assumptions for the exponential Euler method}

\label{subsec:model_assumptions}
Our assumptions will be similar to those in the seminal work
\cite{JK09} on exponential Euler integrators.
%
%
%
\begin{assumption}\label{assum:oper}
There exists a strictly increasing sequence $(\lambda_n)_{n=1}^{\infty}$
of positive real numbers such that $Ae_n=-\lambda_n\,e_n$ for
$n\in{\mathbb N}$ and the linear operator $A:D(A) \rightarrow H$
is given as
\begin{equation*}
Av = -\sum^{\infty}_{n=1} \lambda_n \langle e_n,v\rangle
\,e_n\quad\forall\,v\in D(A)
\end{equation*}
where
\begin{equation*}
D(A) = \left \{v \in H \, \Big| \,  \sum^{\infty}_{n=1}\lambda_n^2\,
|\langle v,e_n\rangle|^2 < \infty \right\}.\quad
\end{equation*}
%
%
\end{assumption}
%
%
%
%
We define the family of interpolation spaces of the operator $A$
for $r \ge 0$ as follows
%
%
\begin{equation}\label{eq:interpolationSpace}
H_r:= D((-A)^r) =\left\{v\in H \, \Big| \,  \sum^{\infty}_{n=1}\lambda_n^{2r}
\,|\langle v,e_n\rangle|^2 < \infty \right\}.
\end{equation}
%
%
\par
The $Q$-Wiener process is defined by
\begin{equation}\label{eq:qwiener}
W_t := \sum^{\infty}_{n=1}\sqrt{q_n}e_n w^n_t,
\end{equation}
where $(w^n_t)_{n=1}^{\infty}$ is a sequence of independent
scalar-valued Wiener processes, and the non-negative sequence
$(q_n )_n \subset [0,\infty)$ satisfies the following:
\begin{assumption}
\label{assum:noise} There exists a constant $\phi \in (0,1)$ such that
\[
\sum^{\infty}_{n=1} (\lambda_n)^{2\phi-1}q_n < \infty.
\]
\end{assumption}
Let $L(H_{r_1},H_{r_2})$ denote the set of bounded linear operators mapping from 
$H_{r_1}$ to $H_{r_2}$, and the let $L(H_r) := L(H_r,H_r)$.
\begin{assumption}
\label{assum:func} The reaction term $f:H \rightarrow H$ is twice
continuously \blue{Fr{\'e}chet} differentiable, where its derivatives satisfy the following
\[
\|f'(x) - f'(y)\|_{L(H)} \leq C\|x-y\|_H, \quad \|(-A)^{-r}f'(x)(-A)^rv\|_H \leq C \|v\|_H,
\]
for all $x,y \in H, v \in D((-A)^{\phi})$, and $r = \{0,1/2,1\}$,
and
\[
\|A^{-1}f''(x)(v,w) \|_H\leq C\|(-A)^{-1/2}v\|_H \|(-A)^{-1/2}w\|_H,
\]
for all $v,w \in H$, where $C>0$ is a positive constant.
\end{assumption}

\begin{assumption}
\label{assum:ic} The initial value $u_0$ is a $D((-A)^{\phi})$-valued random
variable, that satisfies
$$
\mathbb{E}[ \|(-A)^{\phi}u_0\|^4_H ] < \infty,
$$
for the constant $\phi>0$ in Assumption \ref{assum:noise}.
\end{assumption}


\section{Model assumptions for the Milstein method}
\label{subsec:model-assumptions-milstein}
In this section, we present the assumptions for the Milstein method~\cite{jentzenRockner2015} in the setting that is relevant for this paper: when 
the operators $A$ and $Q$ share eigenspace and for reaction-diffusion 
SPDE~\eqref{eq:SPDE} with additive noise. 

\begin{assumption}[Drift coefficient and noise assumption]\label{assum:milstein}
Let Assumption~\ref{assum:oper} hold 
and let Assumption~\ref{assum:noise} hold for some $\phi \in (1/2,1)$.
Let $\kappa \in [0, \phi)$ and let $f:H_{\kappa} \to H$ 
be a twice continuously Fr{\'e}chet differentiable mapping with 
\[
\sup_{x \in H_{\kappa}} \max(\|f'(x)\|_{L(H)}, \, \|f''(x) \|_{L(H_{\kappa} \times H_{\kappa}, H)} ) < \infty.  
\] 
And for a value $\theta \in [\max(\kappa,\phi-1/2) , \phi)$, it holds that 
\[
  \E{ \|(-A)^{\theta} u_0\|_H } < \infty.
\] 
\end{assumption}  

\begin{remark}
Since $H_\kappa$ is dense in $H$, the operator $f'(x) \in L(H_\kappa, H)$ has a unique extension $\tilde f'(x) \in L(H,H)$ and one should 
interpret the operator norm on the extended domain as follows: 
\[
  \|f'(x)\|_{L(H)} := \sup_{ v \in H_\kappa \setminus\{0\} } \frac{ \|f'(x)(v)\|_{H}}{\|v\|_{H}} = \|\tilde f'(x) \|_{L(H)}.
\]
\end{remark}

\begin{remark}\label{rem:notationMilstein} 
The cryptic parameter $\theta$ is an adaptation of~\cite[Assumption 3]{jentzenRockner2015} to the additive-noise setting with $B(u) = I$ and $U_0 = Q^{1/2}(H)$.  
And, working with Hilbert--Schmidt operator norms, \cite[equation (21)]{jentzenRockner2015} is then fulfilled by 
\[
  \|B(u)\|_{HS(U_0, H_{\phi-1/2} )} = \|B(u)Q^{1/2}\|_{HS(H, H_{\phi-1/2} )} := 
  \sum_{n=1}^\infty q_n \lambda_n^{2 \bar \phi-1} < \infty. 
\] 
What we represent by $\kappa$, $\phi -1/2$ and $\theta$ is respectively denoted by $\beta$, $\delta$ and $\gamma$ in~\cite{jentzenRockner2015}. 
\cite[equation (22)]{jentzenRockner2015} is trivially fulfilled 
since $B'(u) =0$ and choosing, in the paper's notation, $\alpha =0$ and 
$\vartheta = \max(1/2 - \phi, 1/4)$, it follows that~\cite[equation (23)]{jentzenRockner2015} 
holds for any $\theta \in [ \max(\kappa,\phi-1/2), \phi)$, since
\[
  \begin{split}
  \|(-A)^{-\vartheta} B(u) Q^{-\alpha} \|_{HS(U_0, H)} &
  = \|(-A)^{-\max(1/2-\phi,1/4)}Q^{1/2} \|_{HS(H, H)}\\
  & \le \sum_{n=1}^{\infty} \lambda_{n}^{2\phi-1} q_n \stackrel{Assumpt.~\ref{assum:noise}}{<} \infty.
  \end{split}  
\]
\cite[equation (23)]{jentzenRockner2015} does indeed not depend on the value $\theta$ in the additive-noise setting, but $\theta$ does enter as a constraint on the regularity of the initial data in~\cite[Assumption 4]{jentzenRockner2015}.
Our lower bound $\phi > 1/2$ is due to the constraint $\delta >0$
in~\cite[Assumption 3]{jentzenRockner2015} and our upper bound $\kappa < \phi$
is due to the constraint $\beta < \delta +1/2$ in~\cite[Assumption 3]{jentzenRockner2015}. 
\end{remark}

\section*{Acknowledgments}
Research reported in this publication received support from the
Alexander von Humboldt Foundation. NKC and AJ are sponsored by KAUST
baseline funding, HH acknowledges support by RWTH Aachen University.

\end{document}